\documentclass[11 pt, reqno]{amsart}

\usepackage{microtype}
\usepackage{amsthm}
\usepackage{amsmath}
\usepackage{amssymb}
\usepackage{amsfonts}
\usepackage{latexsym}
\usepackage{hyperref}
\usepackage{bbm}
\usepackage{comment} 
\usepackage{enumitem}
\usepackage{cite}
\usepackage{braket}

\usepackage{chngcntr}
\counterwithin*{equation}{section}
\usepackage{mathtools}
\newcommand{\defeq }{\vcentcolon=}

\usepackage[letterpaper, left=1.2in, right=1.2in, top=1.2in, bottom=1.2in]{geometry}
\usepackage{color}
\usepackage{fancyhdr}
\usepackage{latexsym}

\newtheorem{ccounter}{ccounter}[section]
\newtheorem{thm}[ccounter]{Theorem}
\newtheorem{lem}[ccounter]{Lemma}
\newtheorem{cor}[ccounter]{Corollary}
\newtheorem{defn}[ccounter]{Definition}
\newtheorem{prop}[ccounter]{Proposition}
\newtheorem{ass}[ccounter]{Assumption}
\newtheorem{ex}[ccounter]{Example}

\theoremstyle{remark}

\newtheorem{remark}[ccounter]{Remark}

\def\bet{\begin{thm}}
\def\eet{\end{thm}}
\def\bel{\begin{lem}}
\def\eel{\end{lem}}
\def\bas{\begin{ass}}
\def\eas{\end{ass}}
\def\bec{\begin{cor}}
\def\eec{\end{cor}}
\def\bed{\begin{defn}}
\def\eed{\end{defn}}
\def\bep{\begin{prop}}
\def\eep{\end{prop}}
\def\beq{\begin{equation}}
\def\eeq{\end{equation}}
\def\bea{\begin{equation*}}
\def\eea{\end{equation*}}
\def\tr{\mathrm{tr}}
\def\bex{\begin{ex}}
\def\eex{\end{ex}}
\def\ep{\qed \\}
\def\be{\beq}
\def\ee{\eeq}
\def\ber{\color{red}}
\def\eer{\normalcolor}

\def\benr{\begin{enumerate}[label=(\roman*)]}
\def\eenr{\end{enumerate}}

\def\N{\mathbb{N}}
\def\C{\mathbb{C}}

\def\R{\mathbb{R}}
\def\P{\mathbb{P}}
\def\E{\mathbb{E}}

\def\fc{m_{\mathrm{fc}}}
\def\sc{\rho_{\mathrm{sc}} }

\def\one{{\mathbbm 1}}
\def\eps{\varepsilon}

\def\Im{\operatorname{Im}}
\def\Re{\operatorname{Re}}

\newcommand{\dd}{d}
\newcommand{\bma}{\begin{bmatrix}}
\newcommand{\ema}{\end{bmatrix}}
\def\trans{*}
\def\tr{\operatorname{Tr}}

\newcommand{\abs}[1]{\lvert #1 \rvert}
\newcommand{\absa}[1]{\left\lvert #1 \right\rvert}

\newcommand{\norm}[1]{\lVert #1 \rVert}
\newcommand{\norma}[1]{\left\lVert #1 \right\rVert}

\def\diag{\operatorname{diag}}
\def\fa{{\mathfrak a}}
\def\fb{{\mathfrak b}}
\def\fc{{\mathfrak c}}
\def\ft{{\mathfrak t}}
\def\ii{i}
\def\supp{\operatorname{supp}}

\newcommand{\hs}[1]{\norm{#1}_{\mathrm{HS}}}

\newcommand{\OO}[1]{\mathcal{O}\!\left(#1\right)}

\newcommand{\half}{\frac{1}{2}}

\begin{document}
\title{Universality of the least singular value for the sum of random matrices}
\author{Ziliang Che}
\thanks{Z.C. is partially supported by NSF grant DMS-1607871.}
\author{Patrick Lopatto}
\thanks{P.L. is partially supported by the NSF Graduate Research Fellowship Program under Grant DGE-1144152.}
\maketitle

{\setcounter{tocdepth}{1}
\tableofcontents
}

\section{Introduction} \label{s:int}

The problem of effectively bounding the least singular value of a random (and possibly non-square) matrix with independent entries has received tremendous attention from mathematicians and computer scientists \cite{RV, basak2017invertibility, LR12,  rebrova2018coverings, tatarko2018upper, rudelson2008invertibility, rudelson2008littlewood, rudelson2009smallest, tikhomirov2015limit,tikhomirov2016smallest,tao2009inverse, ST02, cipolloni2019edge, cipolloni2019optimal}. Recent investigations have also focused on matrices with dependent entries \cite{C16,bose2018smallest,tikhomirov2017invertibility, livshyts2018smallest,  gotze2015minimal, kopel2018random, adamczak2008smallest} and random tensors \cite{vershynin2019concentration, abbe2015reed, anari2018smoothed, bhaskara2014smoothed}. 

A closely related question is to determine the limiting distribution of the least singular value, suitably rescaled, as the size of the matrix tends to infinity. For square matrices with independent entries, it is known that this distribution does not depend on the entry distributions and is equal to the one obtained from a matrix of i.i.d.\ Gaussian random variables (which may be computed exactly). This phenomenon is known as \emph{universality} of the least singular value and was proved for entry distributions with mean zero and variance one in \cite{leastsv} using ideas from the method of property testing in the study of algorithms. In \cite{che2019universality}, universality of the least singular value for square matrices was studied from a dynamical viewpoint and shown to hold for matrices whose entries may be sparse, weakly correlated, and have unequal variances. We also note that the case of genuinely rectangular matrices was taken up in \cite{feldheim2010universality}.

In this work we prove universality of the least singular value for the matrix
\begin{equation}
M =  R^\trans X T+ U^\trans Y V,
\end{equation}
where $X$ and $Y$ are deterministic matrices and $R,T,U,V$ are independently sampled from the Haar measure on the unitary group. The model $M$ is a natural interpretation of the notion of the sum of two generic random matrices and exhibits strong correlations between its entries, unlike the matrices studied in \cite{leastsv, che2019universality}. It was previously studied in \cite{bao2019local}, where its spectrum was controlled on scale $N^{-1 + \eps}$. 

The Hermitian version of this model, $H = V^* X V + U^* Y U$,  has attracted significant interest. The weak convergence of the empirical distribution was obtained first in \cite{voiculescu1991limit} and later shown in \cite{speicher1993free, pastur2000law, collins2003moments, biane1998representations} using alternative techniques. Convergence was then investigated on scales decaying in $N$ in \cite{kargin2012concentration,kargin2015subordination} and established on the optimal scale $N^{-1 + \eps}$ through the series of works \cite{bao2017convergence, bao2017spectral, bao2017local,bao2016local}. The latter results were used to show universality of local spectral statistics in the bulk of the spectrum \cite{che2017local}.

Our work is based on the dynamical approach to random matrix theory developed in the last decade. Based on resolvent estimates and a precise analysis of the short-time behavior of Dyson Brownian motion \cite{Dyson, landon2019fixed, LY, EY15, bourgade2018extreme}, it has succeeded in its original goal of establishing the universality of local spectral statistics for Wigner matrices \cite{wigfixed, EPR10, erdos2009bulk, EYY12, ESY11, localrelax1, localrelax2, EYbook, lee2014necessary}, and has since been applied to investigate universality for numerous other random matrix models. These include random graph models \cite{AC15,BKY, bauerschmidt2017bulk,EKYY13, EKYY12 , bauerschmidt2016local, HLY15}, general Wigner-type matrices \cite{ajanki2017singularities,ajanki2017universality}, band matrices \cite{bourgade2016universality, bourgade2018random, bourgade2019random, yang2018random}, and matrices with few moments \cite{aggarwal2019bulk,aggarwal2018goe}. 

Our proof follows closely the method used in \cite{che2017local} to show universality for the Hermitian model. Two primary inputs in that work were a carefully chosen flow $U(t)$ on the unitary group which leaves the eigenvalue distribution of $H$ unchanged but produces a system of SDEs for the eigenvalue process closely resembling Dyson Brownian motion, and a weak local law throughout the spectrum (including the spectral edges) which was used as an \emph{a priori} input to study the flow of the eigenvalues. We show how similar inputs may be obtained for the model $M$ through a slightly more involved analysis, which proceeds by transforming the problem from one about the singular values of a $N\times N$ non-Hermitian matrix to the eigenvalues of a $2N \times 2N$ Hermitian matrix. The resulting eigenvalue flow is not a Dyson Brownian motion, but instead similar to a symmetrized version studied in \cite{che2019universality}, and the short-time relaxation result for the symmetrized flow in that work is a crucial input here.

Compared to \cite{che2017local}, we derive the weak local law in a slightly different way, involving a general stability analysis of the system of equations that define the free convolution of two measures. While the essential technical content is unchanged, this somewhat streamlines the proofs. Further, we use \cite{bao2019local} to prove a strong law at small energies, paralleling the use of \cite{bao2017convergence, bao2017spectral, bao2017local,bao2016local} in \cite{che2017local} to establish a strong law in the bulk of the spectrum. We also comment on an interesting difference between the real and complex cases which does not arise in the Hermitian model. \\

\noindent {\bf Acknowledgments.} The authors thank Benjamin Landon for comments on a preliminary draft of this paper. They also thank the anonymous referees for their detailed comments, which substantially improved the paper.

\section{Overview and main result} \label{s:overview}

\subsection{Overview} In this section we define the model under consideration and state our main result. The main technical input is Theorem \ref{t:shorttime} about short-time universality for the singular values of the model as it undergoes a time-dependent perturbation. Its proof occupies the bulk of this work. In Section \ref{s:dynamicsdefinition} we define this perturbation and the associated stochastic differential equations governing the evolution of the singular values. In Sections \ref{s:locallaw} and \ref{s:evbounds} we prove various estimates necessary to study the short-time behavior of these SDEs. Their well posedness and the fact that they represent the claimed singular value evolution are proved in Section \ref{s:wp}. In Section \ref{s:sdeanalysis} we compare the SDEs for the singular values to a symmetrized Dyson Brownian Motion flow. The short-time behavior of this flow was studied in \cite{che2019universality}, and by combining our comparison with the main result of that work, we achieve a proof of Theorem \ref{t:shorttime}. As corollary we deduce our main result, Theorem \ref{t:main}.  Appendix \ref{a:ito} contains a computation using It\^o's formula that is required for Section \ref{s:wp}.

For concreteness, we focus on a model where deterministic initial data is conjugated by unitary matrices. It is natural to also consider the analogous model with conjugation by orthogonal matrices. The SDEs for the evolution of the singular values in the second case lack a certain influential repulsion term compared to the first, and as a result the least singular value displays qualitatively different behavior. (This distinction was already noted in \cite{che2019universality} for a different ensemble). Fortunately, our methods suffice to treat this case too. The difference between the behavior of the least singular value in the real and complex models and the necessary modifications to the proof are discussed in Appendix \ref{a:realcase}. 

Finally, Appendix \ref{s:estimates} contains some preliminary estimates required for our analysis. 

\subsection{Model}\label{s:model}
Define 
\begin{equation}\label{e:model}
M =  R^\trans X T+ U^\trans Y V,
\end{equation}
where $X=\diag(x_1,\dots , x_N)$ and $Y=\diag(y_1, \dots , y_N)$ are deterministic diagonal matrices and $R,T,U,V$ are independent and distributed according to the Haar measure on the unitary group $U(N)$. We suppose that 
\begin{equation}\label{e:xybound}
0\le x_i \le C_0, \quad 0 \le y_i \le C_0
\end{equation}
 for some constant $C_0$ independent of $N$. Denote the empirical measures of $X$ and $Y$ by 
\begin{equation}
	\mu_X = \frac{1}{N}\sum_{1\leq i \leq N}\delta_{x_i},\quad \mu_Y = \frac{1}{N}\sum_{1\leq i \leq N}\delta_{y_i}.
\end{equation}

For integers $1\le i \le N$, we define 
\begin{equation}
x_{-i} = - x_{i}, \quad y_{-i} = - y_i.
\end{equation}
We denote the symmetrized version of the empirical  measures of $X$ and $Y$ by\footnote{For a general Borel measure $\nu$ we define $\nu^{\mathrm{sym} }(A) = \frac{1}{2} [\mu(A) + \mu(-A)] $ for any Borel set $A\subset \mathbb R$. }
\begin{equation}
	\mu^{\mathrm{sym}}_X = \frac{1}{2N}\sum_{1\leq\abs{i}\leq N}\delta_{x_i},\quad \mu^{\mathrm{sym}}_Y = \frac{1}{2N}\sum_{1\leq\abs{i}\leq N}\delta_{y_i}.
\end{equation}
Let $\C^+ = \{ z : \Im z > 0 \}$ denote the complex upper half plane.
For $ z\in \C^+$, define the Stietjes transforms
\begin{equation}\label{e:mxz}
	m_X(z)= \int_\R \frac{\mu^{\mathrm{sym}}_X(\dd x)}{x-z},\quad 
	m_Y(z)= \int_\R \frac{\mu^{\mathrm{sym}}_Y(\dd y)}{y-z}.
\end{equation}
We assume:
\begin{enumerate}
	\item  For any $a>0$, there is a constant $C_a>0$ independent of $N$ such that
		\begin{equation}\label{e:mYbound}
			\sup_{E\in\R, \, \eta\geq N^{-1+a}} \absa{m_Y(E+\ii\eta)} \leq C_a.
		\end{equation}
	\item There are compactly supported probability measures $\mu_1$, $\mu_2$ such that $\mu_X\to \mu_1$ and $\mu_Y\to \mu_2$ weakly, and at least one of $\mu_1, \mu_2$ has a bounded Stieltjes transform.\footnote{Observe this implies $\mu^{\mathrm{sym}}_X\to \mu^{\mathrm{sym}}_1$ and $\mu^{\mathrm{sym}}_Y\to \mu^{\mathrm{sym}}_2$ weakly.}
	
	\item Neither of $\mu^{\mathrm{sym}}_1, \mu^{\mathrm{sym}}_2$ is a single point mass, and at least one is supported at more than two points. 
	
	\item The Stieltjes transform of the measure $\mu_X$ converges to that of $\mu_1$ with polynomial speed, in the sense that there exists a constant $\fc_X>0$ such that 
	\begin{equation}
	\left|  \frac{1}{N} \sum_{i=1}^N  \frac{1}{x_i - E - i\eta} - \int \frac{d \mu_1(x)}{x - E - i \eta}   \right| \le N^{-\fc_X}
	\end{equation}
	for $\eta \ge N^{-\fc_X}$. 
	\item\label{strongassumption} The particle $y_k$ is close to the deterministic location $y^*_k$ in the sense that for any $c>0$,
	\begin{equation}
	\sup_{1\le k \le N } |y_k - y^*_k | \le N^{-1 +c },
	\end{equation}
	where $y^*_k$ is the $k$-th $N$-quantile of $\mu_2$ defined by 
	\begin{equation}
	y^*_k = \inf \left\{ s \colon \int_{-\infty}^s \mu_2 (dy) = \frac{k}{N} \right\}.
	\end{equation}
	\item\label{edges} The measure $\mu_2$ has a continuous density and there are constants $\fc,  \delta_0 >0$ such that for any $x \in \supp \mu_2$ and $0 \le h \le \delta_0$, 
	\begin{equation}
	\mu_2 ( [x - h , x + h ]  ) \ge h^{2- \fc}.
	\end{equation}
	\item\label{densityatzero} The free convolution\footnote{We recall the definition of free convolution in Subsection \ref{s:weaklaw1}} $\mu^{\mathrm{sym}}_1  \boxplus \mu^{\mathrm{sym}}_2$ has a density $\rho(x)$ in a neighborhood of zero such that 
	\begin{equation}\label{e:fix0}
	c \le  \rho(x) \le C, \quad \rho(0) = \frac{1}{\pi}
	\end{equation}
	for some constants $C, c >0$ and all $x \in [ - c , c]$. 
\end{enumerate}
The first assumption is to prevent the $y_i$ from accumulating around any point $E\in\R$.  This is illustrated by Proposition \ref{prop:empiricalbounds}. The second and third are required to apply \cite[Theorem 4.4]{bao2019local} to control the Green functions of (a modification of) $M$, as is done in Section \ref{s:evbounds}. The remaining assumptions are required in Subsection \ref{s:deterministicestimates}.

Assumption \eqref{strongassumption} requires that $Y$ obey a strict rigidity condition. However, it can often be relaxed in practice, for instance near the the spectral edges of $\mu_2$. In particular, $Y$ can be taken to be the spectrum of a Wigner matrix, or more generally a matrix of general Wigner-type, as in \cite{ajanki2017universality}. We refer the reader to the remark in \cite[Subsection 2.1]{che2017local} and \cite[Corollary 2.11]{ajanki2017universality} for more on this point.

The condition \eqref{edges} is technical and says, roughly, that the spectral edges of $\mu_2$ behave sub-linearly, for example like the edges of the semicircle distribution $\sc(x) = (2\pi)^{-1} \sqrt{4-x^2}\, dx$. While this is a strong condition, it is true for a broad class of spectral distributions arising in random matrix theory, including those coming from a matrix of general Wigner-type \cite[Theorem 4.1]{ajanki2017universality}; see \cite[Theorem 2.6]{ajanki2} or \cite[Theorem 2.6]{ajanki2017singularities} for more.

The assumption \eqref{densityatzero} is difficult to check in general. For example, the case where $\mu_1$ is a point mass was considered in \cite[Theorem 2.2]{bao2019local}, whose proof is quite technical. In Appendix \ref{s:estimates} we prove two simple sufficient conditions for \eqref{densityatzero}: both $\mu^{\mathrm{sym}}_1$ and $\mu^{\mathrm{sym}}_2$ have positive density at zero, or $\mu^{\mathrm{sym}}_1=\mu^{\mathrm{sym}}_2$. 

The second equation in \eqref{e:fix0} is necessary to make the scale of the smallest singular value match that of the analogous Gaussian ensemble. We include it for technical convenience, but it could be trivially removed by an appropriate rescaling of $X$ and $Y$. We use it in Section \ref{s:sdeanalysis}.

\subsection{Main result} The following is our main result. It is proved at the end of Section \ref{s:sdeanalysis}.

\bet\label {t:main}
Let $\lambda_1(M_N)$ be the least singular value of the random matrix ensemble \eqref{e:model} defined in Section \ref{s:model}. For all $r\ge 0$, we have 
\beq \P (N \lambda_1 ( M_N) \le r ) = 1 - e^{-r^2} + O(N^{-c}),
\eeq
where $c>0$ is an absolute constant uniform in $r$.
\eet

\section{Definition of dynamics}\label{s:dynamicsdefinition}

\subsection{Unitary Brownian motion}
We use the following definitions. Recall a standard complex Gaussian random variable is such that its real and imaginary parts are independent mean zero normal distributions with variance $1/2$.  
\bed\label{d:cbm}
A complex-valued stochastic process $B(t)$ is called a \emph{standard complex Brownian motion} if $(\sqrt{2}\Im B(t), \sqrt{2}\Re B(t) )$ are independent real standard Brownian motions.
\eed

\bed\label{hbm}
A complex-valued stochastic process $(B_{ij}(t))_{1\leq i,j\leq N}$ is called a \emph{Hermitian Brownian motion} if $(\sqrt{2}\Im B_{ij}, \sqrt{2}\Re B_{ij})_{i<j}$, $(B_{ii})_{1\leq i\leq N}$ are a collection of independent real standard Brownian motions, and $\overline{B}_{ij} = B_{ji}$.
\eed
The following construction parallels \cite[Section 2.3.1]{che2017local}. Given a parameter $\mathfrak{a}\in(0,1)$ we introduce the index set 
\be\label{ifa}
	\mathcal{I}_\mathfrak{a} = \{ (i,j) \colon \absa{y_i-y_j} \geq N^{-1+\mathfrak{a}}\} 
\ee
and let $\mathcal{I}_\mathfrak{a}^c$ be the set of pairs $(i,j)$ not in $\mathcal{I}_\mathfrak{a}$.

We let $U(0) \defeq U$ and $V(0)\defeq V$ evolve according to the following equations. 
\be\label{dU}
	\dd U(t) = \ii  \dd W_1 U(t) -\frac{1}{2}AU(t)\, \dd t,\quad \dd V(t) = \ii \dd W_2 V(t) -\frac{1}{2}  AV(t)\, \dd t
\ee
Here $\dd W_{1}, \dd W_{2},$ and $A$ are defined as follows.  Let $\widetilde W_1$ and $\widetilde W_2$ be independent Hermitian Brownian motions in the sense of Definition \ref{hbm}. For $1\leq i,j\leq N$, define the matrix processes $W_1$ and $W_2$ entrywise by
\begin{align}\label{w12}
	(W_1)_{ij} &=\frac{\one_{(i,j)\in \mathcal{I}_\mathfrak{a}}}{\sqrt{2N}}\left( \frac{1}{\absa{y_i-y_j}} (\widetilde W_1)_{ij} + \frac{1}{y_i+y_j}(\widetilde W_2)_{ij}\right), \\
	(W_2)_{ij} &=\frac{\one_{(i,j)\in \mathcal{I}_\mathfrak{a}}}{\sqrt{2N}}\left(\frac{1}{\absa{y_i-y_j}} (\widetilde W_1)_{ij} -\frac{1}{y_i+y_j}(\widetilde W_2)_{ij}\right).
\end{align}
The diagonal matrix $A$ in \eqref{dU} is given by
\be\label{Aii}
	A_{ii} = \frac{1}{2 N} \sum_{j:(i,j)\in\mathcal{I}_\mathfrak{a}} \left( \frac{1}{(y_i-y_j)^2}+\frac{1}{(y_i+y_j)^2}\right).
\ee

Let us explain why these choices are made. With this definition of $W_1$ and $W_2$, we see that
\begin{equation}
	(\ii \dd W_1 Y-\ii Y\dd W_2)_{ij} = \frac{\ii \one_{(i,j)\in\mathcal{I}_\mathfrak{a}}}{\sqrt{2N}}\left( \mathrm{sgn}(y_j-y_i)( \dd \widetilde W_1)_{ij} + (\dd \widetilde W_2)_{ij}\right).
\end{equation}
We therefore find by the L\'evy criterion that $\sqrt{N} \left((\ii \dd W_1 Y-\ii Y\dd W_2)_{ij}\right)_{(i,j)\in \mathcal{I}_\mathfrak{a}}$ is a family of independent standard complex Brownian motions. In particular, there is no longer a Hermitian symmetry. We write
\begin{equation}\label{e:newbm}
	 (\ii \dd W_1 Y-\ii Y\dd W_2)_{ij} =\frac{\one_{(i,j)\in\mathcal{I}_\mathfrak{a}}}{\sqrt{N}}\dd \widetilde B_{ij},
\end{equation}
where $(\dd \widetilde B_{ij})_{1\leq i,j\leq N}$ is a family of independent standard complex Brownian motions.
We choose $A$ so that the solutions $U(t)$ and $V(t)$ stay on the unitary group. One can verify this by differentiating $(UU^*)(t)$ using It\^o's formula and the above definitions to see that $d(UU^*)(t)$ is constant.

Having defined $U(t)$ and $V(t)$, we can differentiate $M(t)$ (defined in the obvious way using $U(t)$ and $V(t)$) and use \eqref{e:newbm} to see
\be\label{dM}
	\dd M(t) = \frac{1}{\sqrt{N}}U^*(t) (\one_{(i,j)\in\mathcal{I}_\mathfrak{a}}\dd \widetilde B_{ij})V(t) + U^*(t) \widehat A V(t)\, \dd t,
\ee
where $\widehat A$ is a diagonal matrix whose entries are given by
\be\label{hatA}
	\widehat A_{ii} = \frac{1}{2N}\sum_{j:(i,j)\in\mathcal{I}_\mathfrak{a}} \left(\frac{1}{y_j-y_i} -\frac{1}{y_i+y_j}\right).
\ee

\subsection{Canceling mesoscopic drift}
Let $\tau =N^{-1+\fb}$ for a parameter $\fb >0$ that will be chosen in the next section. It is hard to use \eqref{dM} as written because after time $\tau$, the contribution from the second term will be order $ N^{-1+\fb}$, which is larger than the order of the microscopic statistic we are interested in. Therefore we introduce an auxiliary matrix 
\begin{equation}
	\widehat{M}(t) =  M(t) + (\tau-t)U^*(0)\widehat AV(0).
\end{equation}
The process $\widehat{M}(t)$ has the property that $M(\tau)=\widehat{M}(\tau)$ and
\begin{align}
	\dd \widehat{M} &=\frac{1}{\sqrt{N}}U^*(t)\left(\one_{(i,j)\in\mathcal{I}_\mathfrak{a}}\, \dd \widetilde B_{ij}\right)V(t)  + \left(U^*(t)\widehat AV(t)-U^*(0)\widehat AV(0)\right)\, \dd t\\
	&=  \frac{1}{ \sqrt{N}} d \widehat B   + \left(U^*(t)\widehat AV(t)-U^*(0)\widehat AV(0) \right)\, \dd t - \frac{1}{\sqrt{N}}U^*(t)\left(\one_{(i,j)\in\mathcal{I}^c_\mathfrak{a}}\, \dd \widetilde B_{ij}\right)V(t).
\end{align}
 Here $\widehat B$ is a matrix of standard complex Brownian motions. We show in Section \ref{s:sdeanalysis} that the second term, when integrated from $0$ to $\tau$, is $o(N^{-1})$. This is small enough not to disturb the microscopic scale $O(N^{-1})$.

Formally applying It\^{o}'s formula (see Appendix \ref{a:ito} for details) suggests that the evolution of the eigenvalues of $\bma 0&\widehat{M} \\ \widehat{M}^* & 0 \ema$ is governed by the following system of SDEs: 
\begin{equation}\label{e:mhatsv}
d \lambda_i  = \frac{1}{\sqrt{2N}} dB_i   +  \frac{1}{2 N }\sum_{j \neq i } \frac{1 - \gamma_{ij} }{\lambda_i - \lambda_j }\, dt   + R_i,
\end{equation}
where 
\begin{equation}
R_i =  \Re\left\langle j_i, \left(U^*(t)\widehat AV(t)-U^*(0)\widehat AV(0)\right) k_i \right\rangle\, \dd t +  \frac{1}{\sqrt{N}} \Re\left\langle j_i , \left(U^*(t)\left(\one_{(i,j)\in\mathcal{I}_\mathfrak{a}}\dd \widetilde B_{ij}\right)V(t)\right) k_i \right\rangle,
\end{equation}
\begin{equation}
\gamma_{ij} =  \frac{1}{2} \sum_{(a,b) \in \mathcal{I}_\mathfrak{a}^c} |w_i(a)|^2 |z_j(b)|^2    +  \sum_{(a,b) \in \mathcal{I}_\mathfrak{a}^c}  |w_j(a)|^2 |z_i(b)|^2, \quad w_i = Uj_i , \quad z_i = Vk_i,
\end{equation}
and for $i < 0$ we set $R_{i} = - R_{-i}$ and $\gamma_{ij} = - \gamma_{-i,j}$. Here $j_i$ and $k_i$ are the columns of the matrices $J$ and $K$ in the singular value decomposition $\widehat{M} = J S K^*$ with $S$ diagonal. 

In Section \ref{s:wp}, we justify this formal calculation, proving that the SDE \eqref{e:mhatsv} is well-posed and its solution is the eigenvalue process for $\widehat{M}(t)$.

\section{Local law} \label{s:locallaw}

In this section we  prove a local law that is used in the next section to obtain global control on the quantities $w_\alpha(i)$, $z_\alpha(i)$, and $\gamma_{\alpha\beta}$.  
We fix constants $\fa, \fb$ such that 
\begin{equation}
0 < \fb < \frac{\fa}{100} <10^{-4},
\end{equation}
and let $\tau = N^{-1 + \fb}$ denote the short time we study. Any constant $C$ without further specification is a universal constant that may depend on $\mathfrak{a}$ but not on $N$.  It may change from line to line, but only finitely many times, so that it remains finite. The norm $\| \cdot \|$ on matrices denotes the operator norm as an operator $\ell^2 \rightarrow \ell^2$.

We now write $U$ for $U(t)$ and $V$ for $V(t)$, and define
\begin{equation}\label{e:mathcalH}
 {\mathcal{H}}(t) = \bma 0 & U R^\trans X T V^\trans  + Y + (\tau-t) U U(0)^\trans\widehat A V(0)V^\trans  \\ \left( U R^\trans X T V^\trans  + Y + (\tau-t)U U(0)^\trans\widehat A V(0)V^\trans \right)^\trans  & 0 \ema,
\end{equation}
and  ${\mathcal{G}} = (  {\mathcal{H}} - z )^{-1}$. 
Note that the $2N$ eigenvalues of $\mathcal H (t) $ are exactly the $N$ singular values of $\widehat M(t)$ for $t \in [ 0 , \tau]$, where each singular value appears with positive and negative sign. (The off-diagonal blocks of $\mathcal H$ come from the matrix $\widehat M(t)$ multiplied by $U$ on the left and $V^\trans$ on the right.) We also define
\begin{equation}
\widehat {\mathcal{H}}(t) = \bma 0 & U R^\trans X T V^\trans  + Y + (\tau-t) \widehat A  \\ \left( U R^\trans X T V^\trans  + Y + (\tau-t) \widehat A  \right)^\trans  & 0 \ema,
\end{equation}
and  $\widehat{\mathcal{G}}(t) = (  \widehat {\mathcal{H}}(t) - z )^{-1}$. We will begin by studying $\widehat {\mathcal{H}}(t)$, since the lack of randomness in the term involving $\widehat A$ makes it more tractable. 
We then relate it to $\mathcal H(t)$ in the proof of Corollary \ref{c:glaw}.

\subsection{Concentration of Green's functions}

The main probabilistic tool in this section is the following concentration result about the Haar measure on the unitary group $U(N)$.  We use the following notation for the Hilbert--Schmidt norm of matrices:
\begin{equation}\label{e:hs1}
	\hs{A}= \sqrt{ \tr(AA^*)}.
\end{equation}
We also recall the equivalent characterization of this norm in terms of the sum of the squares of the matrix entries:
\begin{equation}\label{e:hs2}
\hs{A}^2 = \sum_{i,j} |A_{ij}|^2.
\end{equation}
The next proposition follows from a theorem by Gromov; see \cite[Corollary 4.4.28]{AGZ}.
\bep\label{Gromov}
Let $g\colon U(N)\to \C$ be a Lipschitz function with Lipschitz constant $L$ in the Hilbert--Schmidt norm:
\begin{equation}
	\absa{g(X)-g(Y)} \leq L \hs{X-Y},\quad \forall X,Y\in U(N).
\end{equation}
 Let $\P$ be the normalized Haar measure on the unitary group  $U(N)$ and $\E$ be the corresponding expectation. Then there is a constant $c>0$ not depending on $N$ such that
\begin{equation}
	\P[{\absa{g-\E g}}\geq \delta LN^{-\frac{1}{2}}] \leq \exp\left(-{c\delta^2}\right),\quad \forall \delta\geq 0.
\end{equation}
\eep

The above proposition can be applied to the Green's function $\mathcal G$, which is a smooth function of $U$, $V$, and $z$. In particular, the Lipschitz constant of $\mathcal G$ with respect to the variable $U$ or $V$ can be bounded using the imaginary part of $z$, as illustrated by the following propositions.

\bep\label{imG}
	For any $z=E+\ii \eta\in\C_+$ and $1\leq i\leq 2N$, 
	\begin{equation}
		\Im \mathcal G_{ii} \geq \frac{\eta}{(C+\abs{z})^2},\quad \norm{\mathcal G}\leq \eta^{-1}.
	\end{equation}
The same bounds hold for $\widehat {\mathcal G}$.
\eep
\begin{proof}
Let $e_i\in\C^N$ be the unit vector with $1$ on the $i$-th coordinate and $0$ otherwise. Then we have
\begin{equation}
	\Im \mathcal G_{ii} =\Im \langle e_i, \mathcal G e_i\rangle =\frac{1}{2\ii} \langle e_i, (\mathcal G-\mathcal G^\trans )e_i\rangle.
\end{equation}
Note that $\mathcal G-\mathcal G^\trans = \mathcal G^*(z-\bar z)\mathcal G = 2\ii\eta \mathcal G^* \mathcal G$. Therefore, 
\be\label{ward}
	\Im \mathcal G_{ii} = \langle \mathcal Ge_i, \mathcal Ge_i\rangle \eta = \norm{\mathcal Ge_i}^2\eta.
\ee
By the definition of $\mathcal G$ we have  $\norm{e_i}\leq \norm{\mathcal H-z}\norm{\mathcal Ge_i}$. Note that our assumptions and the bound on $\widehat A$ given in \eqref{e:hatAbound} imply that $\norm{\mathcal H-z}\leq C +\abs{z}$.  Then
\begin{equation}
	\Im \mathcal G_{ii} = \norm{\mathcal Ge_i}^2\eta\geq \frac{\eta}{(C+\abs{z})^2}.
\end{equation}
This proves the first inequality. The second inequality follows from the spectral theorem.\end{proof}

\bep\label{Lip}
Fix $z=E+\ii \eta\in\C^+$. Then, for any $1\leq i,j\leq 2N$,  $\mathcal G_{ij}$ is a Lipschitz function of $U$ (or $V$) in the Hilbert--Schmidt norm with Lipschitz constant bounded by $C\eta^{-2}$ for some $C >0$ independent of $N$. The same holds for $\widehat{\mathcal G}_{ij}$.
\eep

\begin{proof}
Let $\widetilde {\mathcal G}$ be the Green's function $\mathcal G$ after replacing $U$ with $\widetilde U\in U(N)$. Then the resolvent identity yields
\begin{equation}
	{\widetilde {\mathcal G}- \mathcal G}  =
	 \mathcal G\bma 0 & (U-\widetilde U) R^\trans X TV^\trans  \\
	   \left( (U-\widetilde U) R^\trans X TV^\trans \right)^\trans& 0\ema \widetilde{ \mathcal G}.
\end{equation}
Therefore, using the general inequalities $\hs{AB} \le \hs{A} \| B\|$ and $\hs{AB} \le \| A \| \hs{B}$,
\be\label{resolvenths}
	\hs{\widetilde {\mathcal G}- \mathcal G}  \leq 2\norm{\mathcal G}\norm{\widetilde {\mathcal G}} \norm{X} \hs{U-\widetilde U}.
\ee
By Proposition \ref{imG}, $\norm{\mathcal G}\leq \eta^{-1}$, and similarly the spectral theorem yields $\norm{\widetilde {\mathcal G}}\leq \eta^{-1}$. By assumption \eqref{e:xybound}  we have $\norma{X}\leq C$. Hence \eqref{resolvenths} implies that
\begin{equation}
	\hs{ {\mathcal G}- \widetilde{\mathcal G} }\leq \frac{2C }{\eta^2}\hs{U-\widetilde U},
\end{equation}
and the conclusion follows.\end{proof}

\subsection{Invariant identities}

Let $E_{ij}$ be the matrix whose $(i,j)$-th entry is $1$ and all the other entries are $0$:
\begin{equation}
	(E_{ij})_{kl} = \delta_{ik}\delta_{jl}.
\end{equation}
The matrix $E_{ij}$ will be either $N$ by $N$ or $2N$ by $2N$, depending on the context. 

For brevity we set $\overline{Y} = Y + (\tau-t) \widehat A$, and we define
\begin{equation}
	\mathcal H_1 = \bma 0 & U R^*X T V^\trans \\ ( U R^*X T V^\trans )^* & 0 \ema, \quad \mathcal H_2 = \bma 0 &  \overline{Y}  \\ \overline{Y}^* & 0 \ema,
\end{equation}
so that $\widehat{\mathcal H}= \mathcal H_1+ \mathcal H_2$. We require the following lemma. 
\bel\label{l:rotinv} For $1\leq i,j\leq N$ or $N+1\leq i,j\leq 2N$, we have
\begin{equation}
	\E\left[ \widehat{\mathcal G} \mathcal H_1 E_{ij} \widehat{\mathcal G}\right] = \E\left[\widehat{\mathcal G}E_{ij} \mathcal H_1 \widehat{\mathcal G} \right].
\end{equation}
\eel
\begin{proof}
 For any $\zeta\in\C$, define an $N \times N$ unitary matrix $Q(\zeta)$ by
\begin{equation}
	Q(\zeta)= \exp(\zeta E_{ij} -\bar\zeta E_{ji}).
\end{equation}
The derivatives of $Q(\zeta)$ and $Q^*(\zeta)$ with respect to $\zeta$ at $\zeta=0$ are 
\begin{equation}
	\partial_\zeta Q(\zeta) |_{\zeta=0}= E_{ij} ,\quad \partial_\zeta Q^*(\zeta) |_{\zeta=0}= -E_{ij}.
\end{equation}
Let $M(\zeta)= Q(\zeta) UR^*X T V^*  +  \overline{Y}$, let $\mathcal H(\zeta)$ be its symmetrization, and let $G(\zeta, z)$ be the Green's function of $\mathcal H(\zeta)$.  By definition, $M(0)=M$, $\mathcal H(0)=\widehat{\mathcal H}$, and $\mathcal G(0,z)=\widehat{\mathcal G}(z)$.  We differentiate $G(\zeta,z)$ and evaluate the derivative at $\zeta=0$ to obtain
\begin{equation}
	\partial_\zeta  G(\zeta,z) |_{\zeta =0} =  \widehat{\mathcal G}\bma 0 & E_{ij} U R^\trans X T V^\trans \\ -( U R^\trans X T V^\trans )^\trans E_{ij} & 0 \ema \widehat{ \mathcal G}.
\end{equation}

Note that the distribution of $Q(\zeta)U$ is invariant with respect to $\zeta$, so  ${\E [\partial_\zeta G(\zeta,z)]=0}$.  Therefore, the above equality yields
\begin{equation}
	\E \left[  \widehat{\mathcal G}\bma 0 &E_{ij} U R^\trans X T V^\trans  \\ -( U R^\trans X T V^\trans )^\trans E_{ij} & 0 \ema \widehat{\mathcal G} \right] =0.
\end{equation}
This can also be rearranged as 
\begin{equation}
	\E\left[\widehat{\mathcal G}E_{i, j} \mathcal H_1 \widehat{\mathcal G} \right] = \E\left[ \widehat{ \mathcal G} \mathcal H_1 E_{i,j}\widehat{ \mathcal G}\right] .
\end{equation}
This proves the case where $1\leq i,j\leq N$. The other case follows from a similar argument after multiplying $Q(\zeta)$ on the right. 
\end{proof}

\subsection{Asymptotic equations}

Let $m_N(z)$ be the Sieltjes transform of $\widehat{\mathcal H}$. In this subsection we provide a system of equations that $m_N(z)$ satisfies asymptotically.  We require the following high probability notation.

\bed\label{calO}
Given two sequences of random variables $(X_N)$ and $(Y_N)$, we write $X=\OO{Y}$ if there are $c_1,c_2,c_3,N_0>0$ that do not depend on $N$ such that for $N>N_0$,
\begin{equation}
	\P[\absa{X_N}\geq c_1Y_N] \leq \exp\left( - c_2 N^{c_3}\right).
\end{equation}
In general, for some index set $\mathcal A$ (possibly $N$-dependent) and families of random variables $(X(\alpha,N))$ and $(Y(\alpha,N))$ with parameters $\alpha\in\mathcal{A}$ and $N\in\N$, we say that $
	X = \OO{Y}$
uniformly in $\alpha$, if there are constants $c_1,c_2,c_3,N_0>0$ that do not depend on $N$ such that for $N>N_0$,
\begin{equation}
	\P[ \exists \alpha\in\mathcal{A}: \absa{X(\alpha,N)}\geq  c_1Y (\alpha,N)] \leq \exp\left( - c_2 N^{c_3}\right).
\end{equation}
\eed

Define 
\begin{equation} \label{wX}
 w_X = - \frac{ \tr(\widehat{\mathcal G }\mathcal H_1)}{\tr \widehat{\mathcal G}}, \quad w_{\overline{Y}} = -\frac{\tr (\widehat{\mathcal G} \mathcal H_2)}{\tr \widehat{ \mathcal G}},
\end{equation}
which are holomorphic functions on $\C^+$. 

\bel \label{replace} Fix $a>0$. For any $z =  E+\ii \eta\in\C^+$ with $\abs{z}\leq \log N$ and $\eta > N^{-1/10 + a}$,  we have
\begin{equation}
	(\widehat{\mathcal G} \mathcal H_1 )_{ki} = \widehat{\mathcal G}_{ki} \left(-w_X +\OO{\eta^{-5}N^{-\half+a}}\right) + \OO{\eta^{-4}N^{-\half+a}},
\end{equation}
uniformly for $1\leq i, k\leq 2N$.  
\eel

\begin{proof}
By Lemma \ref{l:rotinv} combined with Proposition \ref{Gromov} and Proposition \ref{Lip}, for $1\leq i,j\leq N$ or $N+1\leq i,j\leq 2N$, and any $a>0$, $0\leq k,l\leq 2N$, we have, using the general inequality $\hs{AB} \le \| A \| \hs{B}$ and the symmetry of the matrices in question,	\begin{equation}
		(\widehat{\mathcal G} \mathcal H_1)_{ki} \widehat{\mathcal G}_{jl}= \widehat{\mathcal G}_{ki}( \mathcal H_1 \widehat{\mathcal G})_{jl} +\OO{\eta^{-3}N^{-\half+a}}.
	\end{equation}
We used here our assumptions on $X$ to conclude that $\| \mathcal H_1 \| \le C $.

Take $i,j\in[1,N]$ and  let $l=j$, then sum over $j$:
\be\label{invarianterr}
	( \widehat{\mathcal G} \mathcal H_1)_{ki} \frac{1}{N}\sum_{j=1}^N\widehat{ \mathcal G}_{jj} = \widehat{\mathcal G}_{ki}  \frac{1}{N}\sum_{j=1}^N (\mathcal H_1 \widehat{ \mathcal G})_{jj} +\OO{\eta^{-3}N^{-\half+a}}, \forall 1\leq i\leq N, 1\leq k\leq 2N.
\ee
Note that $U R^*X T V^\trans$ has the same distribution as $(U R^*X T V^\trans)^*$, so  $\bma 0 & I \\ I & 0 \ema \widehat{\mathcal G} \bma 0 & I \\ I & 0 \ema$ has the same probability distribution as $\widehat{\mathcal G}$.  It follows that $\widehat{\mathcal G}_{ii}$ has the same distribution as $\widehat{\mathcal G}_{n+i,n+i}$ for $1\leq i\leq N$.  Therefore,
\begin{equation}
	\E\left[\frac{1}{N}\sum_{j=1}^N \widehat{\mathcal G}_{jj} \right] = \E \left[ \frac{1}{2N} \tr \widehat{\mathcal G} \right].
\end{equation}
Then Proposition \ref{Gromov} and Proposition \ref{Lip} imply that the above identity holds without expectation up to a small error term:
\begin{equation}
	\frac{1}{N}\sum_{j=1}^N \widehat{\mathcal G}_{jj}   = \frac{1}{2N} \tr \widehat{\mathcal G} + \OO{\eta^{-2}N^{-\half+a}}.
\end{equation}	
Similarly, we have
\begin{equation}
	\frac{1}{N}\sum_{j=1}^N (\mathcal H_1 \widehat{\mathcal G})_{jj} =\frac{1}{2N}\tr (\mathcal H_1 \widehat{\mathcal G}) + \OO{\eta^{-2}N^{-\half+a}}.
\end{equation}
Recall the definition \eqref{wX}.  Take the quotient of the above two equations and use Proposition \ref{imG} to get 
\begin{equation}
	\frac{\sum_{j=1}^N (\mathcal H_1 \widehat{\mathcal G})_{jj} }{\sum_{j=1}^N \widehat{\mathcal G}_{jj} } = -w_X + \OO{\eta^{-5}N^{-\half+a}}.
\end{equation}
Here we used the assumed lower bound $\eta > N^{-1/10 + a}$ to ensure the error in the denominator is small. 
We also used H\"{o}lder's inequality for Schatten norms and the second bound in \eqref{imG} (which bounds the largest -- in absolute value -- eigenvalue of $\widehat{\mathcal G}$) to conclude that 
\begin{equation}
\left| \tr (\mathcal H_1 \widehat{\mathcal G}) \right| \le  \| \mathcal H_1\|  \| \widehat{\mathcal G}\|_1 \le C N \eta^{-1}.
\end{equation}
where we used that the Schatten $p$-norm with $p = \infty$ is just the operator norm.

Now we go back to \eqref{invarianterr}, plugging in the above equation to see
\begin{equation}
	(\widehat{\mathcal G} \mathcal H_1)_{ki} = \widehat{\mathcal G}_{ki} \left(-w_X +\OO{\eta^{-5}N^{-\half+a}}\right) + \OO{\eta^{-4}N^{-\half+a}},\quad \forall 1\leq i\leq N, 1\leq k\leq 2N,
\end{equation}
where the extra power of $\eta^{-1}$ in the second error term comes from Proposition \ref{imG} after taking the quotient.
This proves the conclusion for $1\leq i\leq N$. Similarly, we can take $i\in[N+1,2N]$ to obtain the same identity for ${N+1\leq i\leq 2N}, {1\leq k\leq 2N}$. This proves the conclusion for fixed $i$ and $k$. The constants in the $\mathcal{O}$ notation are uniform in $i$ and $k$. \end{proof}

Before proceeding to the next lemma, we prove the following technical proposition.
\bep\label{prop:resolventerror}
Let $A,B,R$ be square matrices of the same dimension such that ${\norm{AR}\leq \delta<1/2}$ and 
\begin{equation}
	A(B+R)=I.
\end{equation}
Then
\begin{equation}
	A = B^{-1} +  O(\delta \norma{A}).
\end{equation}
Here $ O(\cdot)$ is in the sense of operator norm.
\eep
\begin{proof}
We immediately have $
	B= A^{-1} (I-AR).$
Hence 
\begin{equation}
	B^{-1} = (I-AR)^{-1} A =\left( \sum_{k=0}^{\infty} (AR)^k \right)A.
\end{equation}
By the assumption that $\norm{AR}<1/2$, we have $\norma{\sum_{k=0}^{\infty} (AR)^k} \leq 2$. Hence,
\begin{equation}
	\norma{B^{-1}}\leq 2 \norma{A}.
\end{equation}
On the other hand, $A(B+R)=I$ implies
\begin{equation}
	A = B^{-1} - ARB^{-1}.
\end{equation}
Note that $\norma{ARB^{-1}}\leq \delta \norma{B^{-1}} \leq 2\delta\norma{A}$. This gives
\begin{equation}
A = B^{-1} + O(\delta \norma{A}).
\end{equation}
\end{proof}

Now we are ready to prove that the three holomorphic functions $m_N$, $w_X$, and $w_{\overline{Y}}$ approximately satisfy a system of equations for $z\in\C^+$. We recall that $m_X$ was defined in \eqref{e:mxz} as the Stieltjes transform of $\mu^{\mathrm{sym}}_X$, and we let $m_{\overline{Y}}$ denote the Stieltjes transform of $\mu^{\mathrm{sym}}_{\overline{Y}}$.

\bel\label{l:approxsystem}
Fix $a>0$. For any $z\in\C^+$ with $\abs{z}\leq \log N$ and $\eta > N^{-1/10 + a}$, we have
\be\label{sce1}
\left\{
\begin{aligned}
	m_N(z) &= m_X(z+w_{\overline{Y}}) + \OO{\eta^{-6}N^{-\half+a}},\\
	m_N(z) &= m_{\overline{Y}}(z+w_X) + \OO{\eta^{-6}N^{-\half+a}},\\
	\frac{1}{m_N(z)} &= -z - w_X-w_{\overline{Y}},
\end{aligned}\right.
\ee
and for $1\le k \le N$, 
\be
\begin{aligned}
\widehat{\mathcal G}_{kk}&=  \half\left(\frac{1}{\overline{Y}_k-z-w_X} +\frac{1}{-\overline{Y}_k-z -w_X}\right)+ \OO{\eta^{-4}N^{-\half+a}},\\
\widehat{\mathcal G}_{N+k, N + k} &=  \half\left(\frac{1}{\overline{Y}_k-z-w_X} +\frac{1}{-\overline{Y}_k-z -w_X}\right) + \OO{\eta^{-4}N^{-\half+a}}.
\end{aligned}
\ee
\eel

\begin{proof} We start with the following identity, which is equivalent to the definition of $\widehat{\mathcal G}$.
\begin{equation}
	\widehat{ \mathcal G}\bma 0 & \overline{Y} \\ \overline{Y} & 0 \ema + \widehat{\mathcal G }\mathcal H_1 -z \widehat{\mathcal G}= I
\end{equation}
Taking the $(k,k)$-th entry of each of the four blocks, we have
\begin{equation}
	\overline{Y}_k \bma \widehat{\mathcal G}_{k,N+k} &\mathcal  G_{kk} \\ \widehat{\mathcal G}_{N+k,N+k} & \widehat{\mathcal G}_{N+k,k}\ema  + \bma (\widehat{\mathcal G} \mathcal H_1)_{kk} & (\widehat{\mathcal G} \mathcal H_1)_{k,N+k} \\ ( \widehat{\mathcal G} \mathcal H_1)_{N+k,k} & (\widehat{\mathcal G} \mathcal H_1)_{N+k,N+k}\ema - z\bma \widehat{\mathcal G}_{kk} & \widehat{\mathcal G}_{k,N+k} \\ \widehat{\mathcal G}_{N+k,k} & \widehat{\mathcal G}_{N+k,N+k}\ema   = \bma 1 & 0\\ 0 & 1 \ema ,
\end{equation}
or equivalently,
\begin{equation}\label{bp1}
	 \bma \widehat{\mathcal G}_{kk} & \widehat{\mathcal G}_{k,N+k} \\ \widehat{\mathcal G}_{N+k,k} & \widehat{\mathcal G}_{N+k,N+k}\ema \bma -z & \overline{Y}_k \\ \overline{Y}_k &-z\ema  + \bma (\widehat{\mathcal G} \mathcal H_1)_{kk} & (\widehat{\mathcal G} \mathcal H_1)_{k,N+k} \\ (\widehat{\mathcal G} \mathcal H_1)_{N+k,k} & (\widehat{\mathcal G} \mathcal H_1)_{N+k,N+k}\ema = \bma 1 & 0\\ 0 & 1 \ema.
\end{equation}
We apply Corollary \ref{replace} to the matrix involving $\widehat{\mathcal G} \mathcal H_1$ to get
\begin{equation}
	\bma \widehat{\mathcal  G}_{kk} & \widehat{\mathcal  G}_{k,N+k} \\ \widehat{\mathcal  G}_{N+k,k} & \widehat{\mathcal  G}_{N+k,N+k} \ema\left( \bma -z-w_X & \overline{Y}_k \\ \overline{Y}_k &-z-w_X\ema+ \OO{\eta^{-5}N^{-\half+a}}\right) + \OO{\eta^{-4}N^{-\half+a}}=  \bma 1 & 0\\ 0 & 1 \ema ,
\end{equation}
where the $\mathcal O$ notation is used entrywise. Using Proposition \ref{prop:resolventerror} and the fact that $\abs{{\widehat{\mathcal G}}_{ij}}\leq \eta^{-1}$ for all $1\leq i,j\leq 2N$, we have
\begin{equation}
	\bma {\widehat{\mathcal G}}_{kk} & \widehat{\mathcal G}_{k,N+k} \\ \widehat{\mathcal G}_{N+k,k} & \widehat{\mathcal G}_{N+k,N+k} \ema = \bma -z-w_X & \overline{Y}_k \\ \overline{Y}_k &-z-w_X\ema^{-1}\left(I+ \OO{\eta^{-6}N^{-\half+a}}\right).
\end{equation}
This can be written explicitly as 
\begin{multline}
	\bma \widehat{\mathcal G}_{kk} & \widehat{\mathcal G}_{k,N+k} \\ \widehat{\mathcal G}_{N+k,k} & \mathcal  G_{N+k,N+k} \ema  \\ = \bma \half\left(\frac{1}{\overline{Y}_k-z-w_X} +\frac{1}{-\overline{Y}_k-z -w_X}\right)& * \\ *&  \half\left(\frac{1}{\overline{Y}_k-z-w_X} +\frac{1}{-\overline{Y}_k-z -w_X}\right)\ema + \OO{\eta^{-6}N^{-\half+a}},
\end{multline}
where we omit the off-diagonal terms and write them as $*$. We sum over the diagonal terms to get 
\begin{equation}
	\frac{1}{2N}\tr \widehat{\mathcal G} = m_{\overline{Y}} (z+w_X) + \OO{\eta^{-4}N^{-\half+a}}.
\end{equation}
This proves the second equation. To prove the second equation, one replaces $\overline{Y}$ with $RU\overline{Y}V^\trans T^*$ and  $U R^*X T V^\trans$ with $X$ in the definition of $\mathcal H$, then repeats the entire argument word for word. Note that this replacement does not change the definition of $w_{\overline{Y}}$, $w_X$, and $m_N$, since the trace  of a matrix does not change under unitary conjugation.  

Finally, the third equation follows from 
\begin{equation}
	\tr (\widehat{\mathcal G} \mathcal H_1) +\tr(\widehat{\mathcal G} \mathcal H_2) -z\tr \widehat{\mathcal G}=2N,
\end{equation}
which follows from the definition of $\widehat{\mathcal G}$.
\end{proof}

\subsection{Weak law}\label{s:weaklaw1}
In this subsection we prove a weak law for the $\mathcal G_{ii}$.  We use the term \emph{weak} because the result is only valid in the regime $\Im z\geq N^{-c}$ for some small constant $c$.  Thus, it is only slightly stronger than the weak convergence of the corresponding measure $\mu_N$.  Nevertheless, the weak law provides necessary bounds for the eigenvectors of $\mathcal H$.  

We deal with equation \eqref{sce1} in a general setting.  Let $m_\alpha$, $m_\beta$ be the Stieltjes transforms of probability measures $\mu_\alpha$, $\mu_\beta$.  Consider the following deterministic equations for fixed $z\in\C^+$.
\be\label{sce3}
\left\{
\begin{aligned}
	m&= m_\alpha(z+w_\alpha) \\
	m &= m_\beta(z+w_\beta) \\
	\frac{1}{m} &= -z - w_\alpha-w_\beta
\end{aligned}\right.
\ee
Observe that equation \eqref{sce1} is a special case of the above equations plus some error terms. The existence and uniqueness of a solution to this system is known. We call the measure $\mu_\alpha\boxplus\mu_\beta$ given by the following proposition the \emph{free convolution} of $\mu_\alpha$ and $\mu_\beta$. 

\bep\label{p:extension}
Given two probability measures $\mu_\alpha$ and $\mu_\beta$ on $\R$, there exists unique analytic functions $w_\alpha,w_\beta,m\colon\C^+\to\C^+$ satisfying \eqref{sce3}, where $m$ is the Stietljes transform of a probability measure we denote $\mu_\alpha\boxplus\mu_\beta$. Further, suppose $\mu_\alpha\boxplus\mu_\beta$ has a density on an interval $I \subset \R$ that is bounded away from zero. If $\mu_\alpha,\mu_\beta$ are compactly supported,  and none of them is a point mass, then $w_\alpha,w_\beta$ extend continuously to $I$; in particular, $|w_\alpha| \vee |w_\beta|$ is uniformly bounded on compact subsets of $\mathbb{C}^+ \cup I$. If in addition, $\mu_\alpha(\{a\})+\mu_\beta(\{b\})<1$ for all $a,b\in\R$, then $\mu_\alpha\boxplus\mu_\beta$ has a continuous density on $\R$.
\eep
\begin{proof}
    For existence and uniqueness see \cite[Theorem 4.1]{BB07}. For the continuous extension see \cite[Remark 2.4]{Be06} or \cite[Theorem 3.3]{belinschi2008lebesgue}. For the continuous density claim see \cite[Corollary 8]{Be14}.
    
    We remark that the referenced works permit the continuous boundary extensions of $w_\alpha$ and $w_\beta$ to take the value $\infty$, but our hypothesis on the density of $\mu_\alpha\boxplus\mu_\beta$ on $I$ rules this out. Indeed, \cite[Theorem 2.3]{Be06} shows that $m$ also extends continuously to $I$ (again with values that may be $\infty$), and this hypothesis shows that $\Im m(a) > c$ uniformly for $a \in I$ and some constant $c>0$. The uniform boundedness claim then follows from the last equation of \eqref{sce3}.
\end{proof}

We need the stability of the solution to \eqref{sce3} under perturbation. To investigate this, it is convenient to write the equation in a more symmetric form.  In fact, the above equation can be rephrased in terms of $w_\alpha$ and $w_\beta$ only. Define
\begin{equation}\label{e:thehatdefs}
	\hat m_\alpha(\zeta)=-\zeta-\frac{1}{m_\alpha(\zeta)},\quad \hat m_\beta(\zeta)=-\zeta-\frac{1}{m_\beta(\zeta)}.
\end{equation}
Proposition 2.2 in \cite{Ma92} says that $\hat m_\alpha$ and $\hat m_\beta$ are Stieltjes transforms of Borel measures  $\hat{\mu}_\alpha$ and $\hat\mu_\beta$ on $\R$, whose total masses are $\sigma_\alpha^2\defeq \int t^2\mu_\alpha(\dd t)$ and $\sigma_\beta^2\defeq \int t^2 \mu_\beta(\dd t)$ respectively (which are $\leq C$ by assumption \eqref{e:xybound}). 

After eliminating $m$, equation \eqref{sce3} becomes
\be\label{sce4}
\left\{
\begin{aligned}
	w_\alpha &= \hat m_\beta(z+w_\beta)\\
	w_\beta &=  \hat m_\alpha(z+w_\alpha).
\end{aligned}\right.
\ee
Define a function $\Phi\colon (\C^+)^2\to(\C^+)^2$ by
\begin{equation}\label{e:phidef}
	\Phi(\zeta_1,\zeta_2)= (\zeta_1 - \hat m_\beta(z+\zeta_2), \zeta_2 - \hat m_\alpha(z+\zeta_1)).
\end{equation}
Note $\Phi$ depends on a choice of $z \in \C^+$, but we omit this in the notation. The equation \eqref{sce4} is equivalent to
\be\label{sce5}
	\Phi(w_\alpha,w_\beta)=(0,0).
\ee
To show stability, we want to show that the solution to the following perturbed version is close to the solution of \eqref{sce5}, when the perturbation $(r_1,r_2)$ is small enough:
\be\label{sce6}
	\Phi(w_\alpha',w_\beta')=(r_1,r_2).
\ee
It suffices to prove that the matrix of first derivatives $D \Phi$ is non-degenerate at $(w_\alpha,w_\beta)$, i.e.,  $(D\Phi)^{-1}$ is bounded.
\bep\label{p:phiinverse} For all $z\in \mathbb \C^+$,
\be\label{inversedphi}
	\norma{\left(D\Phi(w_\alpha,w_\beta)\right)^{-1}} \leq C (1\vee\eta^{-4}).
\ee
\eep
\begin{proof}
Note that
\begin{equation}\label{e:usethis}
	D\Phi(\zeta_1,\zeta_2) = \bma 1 & -\int_\R \frac{\hat\mu_\beta(\dd x)}{(x-z-\zeta_2)^2} \\ -\int_\R \frac{\hat\mu_\alpha(\dd x)}{(x-z-\zeta_1)^2}  & 1\ema.
\end{equation}
To simplify notation, we define
\begin{equation}
	p= \int_\R \frac{\hat\mu_\beta(\dd x)}{(x-z-w_\beta)^2},\quad q= \int_\R \frac{\hat\mu_\alpha(\dd x)}{(x-z-w_\alpha)^2} ,
\end{equation}
and 
\begin{equation}
	\tilde p = \int_\R \frac{\hat\mu_\beta(\dd x)}{\abs{x-z-w_\beta}^2} ,\quad \tilde q  = \int_\R \frac{\hat\mu_\alpha(\dd x)}{\abs{x-z-w_\alpha}^2}.
\end{equation}
Note $\abs{p}\leq \tilde p\leq \eta^{-2} ,\abs{q}\leq \tilde{q}\leq \eta^{-2}$ and 
\begin{equation}
	D\Phi(w_\alpha,w_\beta) = \bma 1 & -p \\ -q  & 1\ema.
\end{equation}
Therefore we can take the inverse of $D\Phi$ explicitly:
\begin{equation}\label{e:back1}
	\left(D\Phi(w_\alpha,w_\beta)\right)^{-1} = \frac{1}{1-pq} \bma 1 & p \\ q & 1\ema.
\end{equation}
Hence by the elementary bound $\norm{A}\leq \hs{A}$ and the triangle inequality,
\begin{equation}\label{e:back2}
	\norma{\left(D\Phi(w_\alpha,w_\beta)\right)^{-1}} \leq \frac{\sqrt{2+\tilde p ^2+ \tilde q^2}}{1-\tilde p\tilde q}\leq \sqrt{2} \frac{1+\eta^{-2}}{1-\tilde p\tilde q}.
\end{equation}
Taking the imaginary part of equation \eqref{sce4}, we have
\be\label{sce5im}
\left\{
\begin{aligned}
	\Im w_\alpha &= (\eta+\Im w_\beta) \tilde p\\
	\Im w_\beta &= (\eta+\Im w_\alpha)\tilde q.
\end{aligned}\right.
\ee
Therefore,
\begin{equation}
	\tilde p\tilde q = \frac{(\Im w_\alpha)( \Im w_\beta)}{(\eta+\Im w_\alpha)(\eta+\Im w_\beta)} < 1. 
\end{equation}
A simple calculation yields $\frac{1}{1-\tilde p\tilde q} \leq (1+ \frac{\Im w_\alpha \wedge \Im w_\beta}{\eta})$.  Therefore,
\begin{equation}
	\norma{\left(D\Phi(w_\alpha,w_\beta)\right)^{-1}} \leq \sqrt{2} \left(1+ \frac{\Im w_\alpha \wedge \Im w_\beta}{\eta}\right)(1+\eta^{-2}).
\end{equation}
By equation \eqref{sce4} and the fact that $\hat\mu_\alpha$ and $\hat\mu_\beta$ have total masses $\leq C$, we have ${\Im w_\alpha \wedge \Im w_\beta \leq C \eta^{-1}}$.  Therefore
\be
	\norma{\left(D\Phi(w_\alpha,w_\beta)\right)^{-1}} \leq C (1+\eta^{-2})^2\leq C (1\vee\eta^{-4}).
\ee\end{proof}

\bep\label{p:perturbation}
Fix $z=E+\ii \eta \in\C^+$. Let $w_\alpha,w_\beta\in\C^+$ solve \eqref{sce4}. Let $w_1,w_2, r\in \C^+$ be such that
\begin{equation}
	\Phi(w_1,w_2) = (r_1,r_2)=r.
\end{equation}
Define $\delta w=(w_1-w_\alpha,w_2-w_\beta)$. There exists $c>0$ such that if $\norm{\delta w}_2\leq c(\eta^3\wedge \eta^7)$, then
\begin{equation}
	\norm{\delta w}_2 \leq C (1\vee\eta^{-4}) \norm{r}_2.
\end{equation}
\eep

\begin{proof}
Using \eqref{e:usethis}, it is straightforward to see  that  for fixed $z = E + i \eta$, 
\begin{equation}
	\norm{D\Phi(\zeta_1,\zeta_2)}_\infty\leq  C(1\vee\eta^{-2}),\quad \norm{D^2\Phi(\zeta_1,\zeta_2)}_\infty \leq C \eta^{-3} ,
\end{equation}
for any $(\zeta_1,\zeta_2) \in (\C^+)^2$. By Taylor expansion at $(w_\alpha,w_\beta)$, we have
\begin{equation}\label{e:step111}
	\norm{r -D \Phi(w_\alpha,w_\beta)\delta w}_2 \leq \sup_{(\zeta_1',\zeta_2') \in (\C^+)^2} 4\norm{D^2\Phi(\zeta_1',\zeta_2')}_\infty \norm{\delta w}_\infty^2\leq\frac{C}{\eta^3} \norm{\delta w}^2_2.
\end{equation}
Using \eqref{e:step111} and the bound \eqref{inversedphi},  we have
\begin{equation}
	\norm{(D \Phi(w_\alpha,w_\beta))^{-1} r -\delta w}\le \norm{(D \Phi(w_\alpha,w_\beta))^{-1}}  \norm{r -D \Phi(w_\alpha,w_\beta)\delta w}_2 \leq C \eta^{-3}(1\vee\eta^{-4})\norm{\delta w}^2_2.
\end{equation}
Hence,
\begin{equation}
	\norm{\delta w}_2\le \norm{(D \Phi(w_\alpha,w_\beta))^{-1} r}_2 + \norm{(D \Phi(w_\alpha,w_\beta))^{-1} r -\delta w}_2 \le  C (1\vee\eta^{-4}) \norm{r}_2 + C \eta^{-3}(1\vee\eta^{-4})\norm{\delta w}^2_2.
\end{equation}
Using the condition that $\norm{\delta w}\leq c (\eta^3\wedge \eta^7)$ and choosing $c>0$ small enough, the second term on the right side can be absorbed into the left side. Thus,
\begin{equation}
	\norm{\delta w}_2 \leq C(1\vee\eta^{-4}) \norm{r}_2.
\end{equation}
\end{proof}

Let $K> 0$ be a constant smaller than the constant $c$ in the assumption \eqref{e:fix0}. Define the spectral domain
\begin{equation}\label{e:sigmadef}
\Sigma = \left\{ z = E + i \eta \colon \eta \in [N^{-1/100} , 1] , E \in [ -2K, 2K]   \right\},
\end{equation}
and let $\tilde{w}_X$ and $\tilde{w}_{\overline Y}$ solve
\be\label{e:tildew}
\left\{
\begin{aligned}
	\tilde w_X &= \hat m_X(z+\tilde w_{\overline{Y}})\\
	\tilde w_{\overline{Y}} &=  \hat m_{\overline Y}(z+\tilde w_X).
\end{aligned}\right.
\ee

\bec\label{c:wlaw}
For any $z\in \Sigma$ and $t \in [ 0 ,\tau]$, 
\begin{equation}
| w_X - \tilde w_X | \vee | w_{\overline{Y}} - \tilde w_{\overline{Y}} |  = \OO{\eta^{-12}N^{-1/3}}.
\end{equation}
\eec
\begin{proof}
We restrict the following claims to $z\in \Sigma$. Multiplying the first and third equations of \eqref{sce1} gives
\beq
1 = \left( m_X(z+ w_{\overline{Y}}) + \OO{\eta^{-6}N^{-\half+a}} \right) ( -z - w_X - w_{\overline{Y}}).
\eeq
By Proposition \ref{imG}, $\Im m_N(z) \ge c\eta$, so 
\beq
| -z - w_X - w_{\overline{Y}}| = \frac{1}{|m_N(z)|} \le C \eta^{-1},
\eeq
\beq\label{e:divideme}
1 = m_X(z+ w_{\overline{Y}})( -z - w_X - w_{\overline{Y}}) + \OO{\eta^{-7}N^{-\half+a}}.
\eeq
Similarly, using the first equation of \eqref{sce1} (with $a$ small enough, say $a = 1/6$) and Proposition \ref{imG} we have
\begin{equation}
\frac{1}{\left|m_X(z+ w_{\overline{Y}}) \right|}   \le C\eta^{-1}.
\end{equation}
Dividing \eqref{e:divideme} by $m_X(z+ w_{\overline{Y}})$ and rearranging using $\hat m_X(z+ w_{\overline{Y}}) = - z  -  w_{\overline{Y}} - \frac{1}{m_X(z+ w_{\overline{Y}})} $ yields
\beq
\hat m_X (z + w_{\overline{Y}}) = w_X + \OO{\eta^{-8}N^{-\half+a}}.
\eeq
Analogously, 
\beq
\hat m_Y (z + w_X) = w_{\overline{Y}} + \OO{\eta^{-8}N^{-\half+a}}.
\eeq
The claim now follows from choosing $a = 1/6$ and applying Proposition \ref{p:perturbation}. \end{proof}

\bec\label{c:glaw}
There exists a constant $c(\fb) > 0$ such that with probability at least $1 - e^{ - c N^c }$,
\begin{equation}
\sup_{(z,t) \in \Sigma \times [0, \tau] } \left| \mathcal G_{ii}  -   \half\left(\frac{1}{\overline{Y}_i-z-\tilde w_X} +\frac{1}{-\overline{Y}_i-z -\tilde w_X}\right)  \right| \le C \eta^{-14}N^{-1/3},
\end{equation}
where the $i$ in $\overline{Y}_i$ is taken modulo $N$.
\eec
\begin{proof}
By Corollary \ref{c:wlaw}, $|w_X - \tilde w_X| =  \OO{\eta^{-12}N^{-1/3}}$, and by Lemma \ref{l:approxsystem},
\begin{equation}
\widehat{\mathcal G}_{ii }=  \half\left(\frac{1}{\overline{Y}_i-z-w_X} +\frac{1}{-\overline{Y}_i-z -w_X}\right)+ \OO{\eta^{-6}N^{-\half+a}}.
\end{equation}
Because
\begin{equation}
\left| \frac{1}{\overline{Y}_i-z-w_X}  - \frac{1}{\overline{Y}_i-z-\tilde w_X}  \right| \le \frac{|w_X - \tilde w_X| }{\eta^2}, 
\end{equation}
we obtain 
\begin{equation}
\widehat{\mathcal G}_{ii }=  \half\left(\frac{1}{\overline{Y}_i-z-\tilde w_X} +\frac{1}{-\overline{Y}_i-z -\tilde w_X}\right)+ \OO{\eta^{-14}N^{-1/3}}.
\end{equation}
This statement for fixed $z$ and $t$ may be upgraded to the desired statement uniform over $\Sigma \times [0,\tau]$ by a standard stochastic continuity argument, as indicated in the proof of \cite[Theorem 3.16]{che2017local}. Observe that we use Lemma \ref{l:holder} in place of \cite[Theorem 3.3]{che2017local} in that argument. Finally, the estimate may be transferred from $\widehat{\mathcal G}$ to $\mathcal G$ by using the resolvent expansion and Ward identity, again as in the proof of \cite[Theorem 3.16]{che2017local}.
\end{proof}

\section{Eigenvector estimates} \label{s:evbounds}

\subsection{Global eigenvector bounds}

\bet\label{t:globalevbounds} There are constants $p,c(\fb)  >0$ such that the following holds. With probability $1 - e^{- N^c}$,
\begin{equation}
\sup_{0 \le t \le \tau} \max_{1\le \alpha, i \le N} |w_\alpha(i)|^2 + |z_\alpha(i)|^2 \le N^{-1/p}, \quad \sup_{0\le t \le \tau}   \max_{\beta \neq \alpha} \gamma_{\alpha \beta} \le N^{-2/p + \fa}.
\end{equation}
\eet
\begin{proof} Recall the eigenvectors of $\mathcal H$ are of the form $(w_\alpha, z_\alpha)$ and $(- w_\alpha, z_\alpha)$, corresponding to eigenvalues $\lambda_\alpha$ and $-\lambda_\alpha$ respectively. Recall that $\Sigma$ was defined in \eqref{e:sigmadef}. Then for $(z,t) \in \Sigma \times [0,\tau]$ and $1\le \alpha \le N$, 
\beq
\Im \mathcal G_{\alpha \alpha}(z) = \eta \sum_{i =1}^N  \frac{|w_\alpha(i)|^2   }{| \lambda_i - E |^2 + \eta^2} + \eta \sum_{i =1}^N  \frac{ |z_\alpha(i)|^2  }{| - \lambda_i - E |^2 + \eta^2}.
\eeq 
Setting $E = \lambda_i$, we obtain 
\beq
|w_\alpha(i)|^2   \le \eta \Im \mathcal G_{\alpha \alpha} (\lambda_i + i \eta).
\eeq
Denote $z_i = \lambda_i + i \eta$.  By the local law for $\mathcal G_{\alpha\alpha}$, Corollary \ref{c:glaw}, it suffices to lower bound 
\beq |-z_i - \overline{Y}_\alpha - \Re \tilde w_X|\vee \Im \tilde w_X.
\eeq
and 
\beq |-z_i + \overline{Y}_\alpha - \Re \tilde w_X| \vee \Im \tilde w_X.
\eeq
By Corollary \ref{c:lowerbound}, proved in the next subsection, 
\begin{equation}
 |-z_i + \overline{Y}_\alpha - \Re \tilde w_X| \vee \Im \tilde w_X \ge \eta^{1 - \fc /4}, \quad 
|-z_i - \overline{Y}_\alpha - \Re \tilde w_X| \vee \Im \tilde w_X \ge \eta^{1 - \fc /4}.
\end{equation}
Therefore, using Corollary \ref{c:glaw}, we have
\begin{equation}
|w_\alpha(i)|^2   \le \eta^{\fc /4} + C \eta^{-11}N^{-1/3}
\end{equation}
with exponentially high probability. Taking $\eta = N^{-1/q}$ with $q$ large finishes the proof of the bound on $|w_\alpha(i)|^2$.

The bound on $|z_\alpha(i)|^2$ is analogous. Given the bound on $|w_\alpha(i)|^2 + |z_\alpha(i)|^2$, the bound for $\gamma_{\alpha \beta} $ follows from the definition of $\gamma_{\alpha \beta}$ after taking $p$ large enough.
\end{proof}

\subsection{Deterministic estimates}\label{s:deterministicestimates}
Let $m_1, m_2$ be the Stieltjes transforms of $\mu^{\mathrm{sym}}_1, \mu^{\mathrm{sym}}_2$, and let $w_1, w_2$ be the solution to the system
\be
\left\{
\begin{aligned}\label{e:w1def}
	 w_1 &= \hat m_1(z+ w_2)\\
	 w_2 &=  \hat m_2(z+ w_1).
\end{aligned}\right.
\ee
The proof of the following proposition is the same as \cite[Proposition 3.9]{che2017local}.

\bep\label{p:39} There exists $p>0$ such that if $\Im z \ge N^{-1/p}$, then 
\begin{equation}
|m_1(z) - m_X(z) | \le N^{-1/p}, \quad |m_2(z) - m_{\overline{Y}}(z) | \le N^{-1/p}
\end{equation}
for all $z$.
\eep

Given $p>10$, define the spectral domain
\begin{equation}
\Sigma_p = \left\{ z = E + i \eta \colon \eta \in [N^{-1/p^2} , 1] , E \in [ -2K, 2K]   \right\},
\end{equation}

\bec\label{c:w1wX} There exists a universal constant $p>0$ such that for $z \in \Sigma_p$,
\begin{equation}
| w_1 - \tilde w_X | \vee | w_2 - \tilde w_{\overline{Y}} | \le N^{-1/p} \eta^{-p}, \quad |\tilde w_X| \vee |\tilde w_{\overline{Y}}| \le C.
\end{equation}
\eec

\begin{proof} This follows from Proposition \ref{p:perturbation} with 
\begin{equation}
\delta w = \left(\frac{1}{m_1(z +  \tilde w_{\overline{Y}}) } -  \frac{1}{m_X(z + \tilde w_{\overline{Y}}) }   , \frac{1}{m_2(z +  \tilde w_X) } -  \frac{1}{m_{\overline{Y}}(z + \tilde w_X) }   \right).
\end{equation}
We indicate how to bound the first coordinate; the second is analogous. From Proposition \ref{p:39}, 
\begin{equation}
|\delta w_1| \le \frac{N^{-1/p}}{m_1(z +  \tilde w_{\overline{Y}})m_X(z + \tilde w_{\overline{Y}})}.
\end{equation}
For any Stieltjes transform $m(z)$ of a measure $\mu$,
\begin{equation}\label{e:mlowerbound}
\Im m(z) \ge \frac{\Im z }{(|z| + \sup_{x\in\supp \mu} |x| )^2}.
\end{equation}
Recall from Proposition \ref{p:extension} that $\Im \tilde w_X \wedge \Im \tilde w_Y >0$, and $| \hat m_X (z) | \le C \eta^{-1}$ since $\hat \mu_X$ has finite mass. 
Taking $p$ large and using $| \tilde w_X | \vee |\tilde w_{\overline{Y}}|  \le C\eta^{-1}$, which follows from \eqref{e:tildew} and the preceding comment, this shows
\begin{equation}|\delta w_1| \le C N^{-1/p} \eta^{-6} \le C\eta.\end{equation}
In the last inequality we used the hypothesis that $\eta \ge N^{-1/p^2}$.

Finally, the second bound follows from the first and Proposition \ref{p:extension}.
\end{proof}
Define $\tilde m(z)$ to be the Stieltjes transform of $\mu^{\mathrm{sym}}_X  \boxplus \mu^{\mathrm{sym}}_{\overline Y}$. We recall this means that $w_X, w_{\overline Y}, m_X, m_{\overline{Y}}, \tilde m$ satisfy \eqref{sce3}.
\bec\label{c:mXmY} Under the same assumptions as Corollary \ref{c:w1wX},
\begin{equation}\tilde m(z) \le C. \end{equation}
\eec
\begin{proof}
As noted in Section \ref{s:model}, we assume that either $\mu_1$ or $\mu_2$ has a bounded Stieltjes transform. By Proposition \ref{p:39}, for $\eta \ge N^{-1/p}$,
\begin{equation}
|m_X ( z + \tilde w_{\overline{Y}} ) | \wedge |m_Y ( z +  \tilde w_X ) | \le |m_1 ( z + \tilde w_{\overline{Y}} ) | \wedge |m_2 ( z +  \tilde w_X ) | + 2 N^{-1/p} \le C.
\end{equation}
Using the definition of $\tilde m$ (recall \eqref{sce3}), this completes the proof.
\end{proof}
The following corollary is essentially contained in the proof of \cite[Theorem 3.15]{che2017local}. We include it for completeness. 
\bec\label{c:lowerbound}
There exists $N_0>0$ such that for $z\in \Sigma_p$ and $N \ge N_0$,
\begin{equation}
 |-z + \overline{Y}_\alpha + \Re \tilde w_X| \vee \Im \tilde w_X \ge \eta^{1 - \fc /4}, \quad |-z - \overline{Y}_\alpha + \Re \tilde w_X|  \vee \Im \tilde w_X \ge \eta^{1 - \fc /4}.
\end{equation}
\eec 
\begin{proof}
We first show 
\begin{equation}
|-z + \overline{Y}_\alpha + \Re \tilde w_X| \vee \Im \tilde w_X \ge \eta^{1 - \fc /4}.
\end{equation}
We suppose, for the sake of contradiction, that $|-z + \overline{Y}_\alpha + \Re \tilde w_X| \vee \Im \tilde w_X  < \eta^{1 - \fc /4}$. Recall 
\begin{equation}
\tilde w_X = - z  - \tilde w_{\overline{Y}} - \frac{1}{m_X (z+ \tilde w_{\overline{Y}})}.
\end{equation}
Taking the imaginary part and using Corollary \ref{c:mXmY} together with the definition \eqref{sce3} shows
\begin{equation}
\Im \tilde w_X + \Im \tilde w_{\overline{Y}} = -\eta  + \frac{\Im m_X(z+ \tilde w_{\overline{Y}})}{| m_X(z+ \tilde w_{\overline{Y}}) |^2 } \ge - \eta  + c  \Im m_X(z+ \tilde w_{\overline{Y}})
\end{equation}
for some $c>0$. 
Set $I(\alpha, \eta) = \left| \left\{  \beta \colon \left| \overline{Y}_\beta  - \overline{Y}_\alpha  \right| \le \eta   \right\}   \right|$. By definition,
\begin{equation}
\begin{split}
\Im m_X(z + \tilde w_{\overline{Y}}) &= \frac{1}{2N} \sum_{\beta = -N}^{N} \frac{\eta + \Im \tilde w_{\overline{Y}}}{|-z + \overline{Y}_\beta + \Re \tilde w_X|^2 +  |\eta +  \Im \tilde w_X|^2  } \\
&\ge \frac{1}{2N} \sum_{|\overline{Y}_\beta - \overline{Y}_\alpha| \le \eta } \frac{\eta + \Im \tilde w_{\overline{Y}}}{|-z + \overline{Y}_\beta + \Re \tilde w_X|^2 +  |\eta +  \Im \tilde w_X|^2  } \\
&\ge \frac{I(\alpha, \eta)} {2N } \frac{\eta}{4 \eta^{2 - \fc/2}} = \frac{I(\alpha, \eta) \eta^{\fc/2 -1}}{8N}.
\end{split}
\end{equation}
Using model assumptions \eqref{strongassumption} and \eqref{edges}, applying $\eqref{e:hatAbound}$, and recalling $\tau = N^{-1+\fb}$, we find
\begin{equation}
I(\alpha, \eta) \ge \left| \left\{  \beta \colon \left| y_\beta  - y_\alpha  \right| \le 2\eta/3   \right\}   \right| \ge \frac{N}{2} \mu_2 ([y_\alpha - \eta/2, y_\alpha + \eta/2])\ge  c N \eta^{2-\fc}.
\end{equation}
This implies $\Im m_X(z + \tilde w_{\overline{Y}}) \ge c \eta^{1 - \fc /2}$, from which we conclude using the above work that
\begin{equation}
\Im \tilde w_X + \Im \tilde w_{\overline{Y}} \ge 2 \eta^{1 - \fc/3}
\end{equation}
when $N$ is large enough.
By the assumption that $\Im \tilde w_X  <  \eta^{1 - \fc /4}$, this implies $\Im \tilde w_{\overline{Y}} \ge \eta^{1 - \fc/3}$.  From $\tilde w_X = \hat m_X ( z + \tilde w_{\overline{Y}})$ and Corollary \ref{c:w1wX}, we obtain
\begin{align}
\Im \tilde w_X = \int_\R \frac{(\eta + \Im \tilde w_{\overline{Y}}) \hat \mu_X (dx)}{|x - z + \tilde w_{\overline{Y}}|^2} &\ge  \eta^{1 - \fc/3} \int_\R \frac{ \hat \mu_X (dx)}{|x - z + \tilde w_{\overline{Y}}|^2 }\\  &\ge  \eta^{1 - \fc/3} \int_\R \frac{ \hat \mu_X (dx)}{ c (x^2 + 1)  } \ge c \eta^{1 - \fc/3}\Im \hat m_X(i).
\end{align}
By Proposition \ref{p:39}, $\Im \hat m_X(i) \ge c$, which implies $\Im \tilde w_X  \ge  \eta^{1 - \fc/4}$, contradicting our assumption (for large enough $N$).

Finally, an analogous argument shows
$|-z - \overline{Y}_\alpha + \Re \tilde w_X| \vee \Im \tilde w_X \ge \eta^{1 - \fc /4},$
and this completes the proof.
\end{proof}

\subsection{Bulk eigenvector bounds} \label{s:53}

\bet\label{t:bulklaw} Let $I$ be the interval in assumption \eqref{densityatzero}. Fix $\nu >0$ and set
\begin{equation}
\mathcal D_{I} = \left\{  z = E + i \eta \colon E \in I, N^{-1 + \nu} \le \eta \le 1  \right\}.
\end{equation}
Then
\begin{equation}\label{e:wxtildelower}
\inf_{z \in \mathcal D_{I}} \Im \tilde w_X \ge c,
\end{equation}
and it holds with overwhelming probability that 
\begin{equation}
\sup_{0 \le t \le  \tau } \sup_{z\in \mathcal D_I} \max_{1 \le i \le 2N } \left| \mathcal G_{ii} (z,t)  - \frac{ z + \tilde w_X }{  ( y_i  + (\tau-t) \widehat{A}_{ii} )^2- ( z + \tilde w_X)^2} \right| \le \frac{N^\nu}{\sqrt{N\eta}},
\end{equation}
where the indices in $y_i$ and $\widehat{A}_{ii}$ are taken modulo $N$, and we require $N \ge N_0$ for some $N_0$ depending on $\nu$.
\eet
\begin{proof}

By assumption, the empirical measure of $U R^\trans X T V^\trans$ converges to $\mu_1$ weakly. Using $\tau = o(1)$ and $\eqref{e:hatAbound}$, $(\tau-t) \widehat A$ is negligible and $Y + (\tau-t) \widehat A$ converges to $\mu_2$ weakly. Fix a small $\sigma>0$. Then by Theorem 4.4 of \cite{bao2019local}, for any fixed $t$,
\begin{equation}
\sup_{z \in \mathcal D_I} \max_{1\le i \le N} \left| \widehat{\mathcal{G}}_{ii}(z,t)  - \frac{ z + \tilde w_X }{  ( y_i  + (\tau-t) \widehat{A}_{ii} )^2- ( z + \tilde w_X)^2}\right| \le \frac{N^\sigma}{\sqrt{N\eta}} 
\end{equation}
with overwhelming probability for sufficiently large $N$ not depending on $t$. Further, by Lemma A.2 of \cite{bao2019local}, there exists $c > 0$ such that $\inf_{z \in \mathcal D_I} \Im \tilde w_X  \ge c$ for large enough $N$, independent of $t$. This implies the desired claim for $\widehat{\mathcal G}(z,t)$ at any fixed $t$. Observe there is an implicit dependence of $\tilde w_X$ on $t$. 

This estimate may then be transferred to $\mathcal G_{ii}$, using the resolvent identity, and made uniform in $t$, using a standard stochastic continuity argument, as in \cite[Theorem 3.16]{che2017local}. This completes the proof. \end{proof}

\bec\label{e:57} Let $I$ be the interval in assumption \eqref{densityatzero}. Then for any $\nu >0$, there exists $N_0(\nu)>0$ such that the following estimates hold with overwhelming probability for $N \ge N_0$:
\begin{equation}
\sup_{0 \le t \le \tau} \max_{\lambda_\alpha \in I} \max_{1\le i \le N} |w_\alpha(i)| + |z_\alpha(i)| \le \frac{N^\nu}{\sqrt{N}}, \quad \sup_{0\le t \le \tau}  \max_{\lambda_\alpha \in I} \max_{\beta \neq \alpha} \gamma_{\alpha \beta} \le \frac{N^{\fa + \nu}}{N}.
\end{equation}
\eec
\begin{proof}
The proof is the same as \cite[Corollary 3.17]{che2017local} (using both conclusions of Theorem \ref{t:bulklaw}) after observing that the eigenvectors of $\mathcal H$ are of the form $(w_\alpha, z_\alpha)$ and $(- w_\alpha, z_\alpha)$.
\end{proof}

Define $m_3(z)$ as the Stieltjes transform of $\mu_1^{\mathrm{sym}} \boxplus \mu_2^{\mathrm{sym}}$. The next proof follows \cite[Theorem 3.14]{che2017local}.

\bet\label{t:m3} Let $I$ and $\mathcal D_I$ be as in Theorem \ref{t:bulklaw}. There exist constants $q > 0 $ and $N_0 >0$ such that 
\begin{equation}\label{e:firstineq}\sup_{z \in \mathcal D_I}  | w_1 - \tilde w_X| \vee |w_2 - \tilde w_{\overline{Y}}| \le N^{-1/q}, \quad \sup_{z \in \mathcal D_I} | m_3(z) - \tilde m(z)  | \le N^{-1/q}
\end{equation}
for $N \ge N_0$.
\eet 
\begin{proof} Let $q > 0$ be a constant to be determined later. We define $\Sigma_1 \subset \mathcal D_I$ by 
\begin{equation}
\Sigma_1 = \left\{ z \in \mathcal D_I \colon | w_1 - \tilde w_X | \vee | w_2 - \tilde w_{\overline{Y}} | \le N^{ -1/q}   \right\}.
\end{equation}
By Corollary \ref{c:w1wX}, $\Sigma_1$ is nonempty. Because the functions involved in its definition are continuous, it is closed. Therefore to show $\Sigma_1 = \mathcal D_I$, it suffices to show $\Sigma_1$ is open in $\mathcal D_I$. 

Because the density of $\mu_1^{\mathrm{sym}} \boxplus \mu_2^{\mathrm{sym}}$ is bounded above and below by positive constants, there exists $c>0$ such that $c^{-1} \ge \Im m_3  \ge c$. 
By taking imaginary parts in \eqref{e:w1def} and the last equation of \eqref{sce3},  using \eqref{e:mlowerbound}, and recalling that $|w_1 | \vee |w_2|$ is bounded on compact subset of $\C^+ \cup \R$ by Proposition \ref{p:extension}, we have 
\begin{equation}\label{e:asdf}
\Im w_2 \ge c \Im w_1,\quad  \Im w_1  \ge c \Im w_2, \quad \Im w_1 + \Im w_2 \ge c > 0,
\end{equation}
which implies $\Im w_1 \wedge \Im w_2 \ge c > 0$.\footnote{We also need $\eta = O(N^{-c})$ for some $c>0$ for the third inequality in \eqref{e:asdf}, but the first inequality of \eqref{e:firstineq} for the complementary regime has already been proved in Corollary \ref{c:w1wX}.}
 These lower bounds permit the use of Proposition \ref{p:39} to conclude that on $\Sigma_1$, 
\begin{equation}\label{e:concludeapriori}
| m_1 ( z + \tilde w_{\overline{Y}})  -  m_X ( z + \tilde w_{\overline{Y}}) | \vee | m_2 ( z + \tilde w_X)  -  m_{\overline Y} ( z + \tilde w_X) | \le 2 N^{-1/p}.
\end{equation}
Therefore, using the definition \eqref{e:thehatdefs} and the lower bound $ \Im m_X (z)\wedge \Im m_{\overline Y} (z)  \ge c \eta $ for some $c>0$ (which follows from the definition of the Stieltjes transform as the trace of a Green's function, as in Proposition \ref{imG}), we have 
\begin{equation}\label{e:cap2}
| \hat m_1 ( z + \tilde w_{\overline{Y}})  -  \hat m_X ( z + \tilde w_{\overline{Y}}) | \vee | \hat m_2 ( z + \tilde w_X)  -  \hat m_{\overline Y} ( z + \tilde w_X) | \le C N^{-1/p}.
\end{equation}
We now claim that on $\mathcal D_I$ the stability of the system of equations \eqref{e:w1def} is improved, so that the operator $\Phi$ from \eqref{e:phidef} satisfies 
\be\label{inversedphi2}
	\norma{\left(D\Phi(w_1,w_2)\right)^{-1}} \leq C.
\ee
To see this, one can reinspect the proof of Proposition \ref{p:phiinverse} using the bound ${\Im w_1 \wedge \Im w_2 \ge c > 0}$, which implies $\tilde p \vee \tilde q \le C$, and the bound 
\begin{equation}
\sup_{\Im z \ge 0} |  p  q |  < \sup_{\Im z \ge 0 } | \tilde  p \tilde q| \le 1,
\end{equation} 
which holds because $\hat \mu_2$ is not a point mass and ${\Im w_1 \wedge \Im w_2 \ge c > 0}$. Because $p, q$ are continuous, we find $|1 - pq| \ge c > 0$ on $I$. Repeating \eqref{e:back1} and \eqref{e:back2} with these improved bounds proves the claim.

Similar reasoning gives $\| D \Phi (\zeta_1, \zeta_2) \|_\infty < C$, $\| D^2 \Phi (\zeta_1, \zeta_2) \|_\infty \le C$ on $\Sigma_1$. We can therefore  repeat the reasoning of the proof of Proposition \ref{p:perturbation} to show that for any $z\in \Sigma_1$, $| w_1 - \tilde w_X | \vee | w_2 - \tilde w_{\overline{Y}} | \le C N^{ -1/p} $. 
The remainder term $r$ in that proof is now bounded using \eqref{e:cap2}.
Therefore, there is a neighborhood of $z$ such that $| w_1 - \tilde w_X | \vee | w_2 - \tilde w_{\overline{Y}} | \le N^{ -1/q} $ when $q>p$. This shows that $\Sigma_1 = \mathcal D_I$. 

Finally, on $\mathcal D_I$ we have using \eqref{sce3}, the lower bound on $\Im w_2$, and Proposition \ref{p:39}, 
\begin{equation}
| m_3(z) - \tilde m(z) | \le  |m_X( z + \tilde w_{\overline{Y}}) - m_1( z + \tilde w_{\overline{Y}})  |    + |m_1( z + \tilde w_{\overline{Y}}) - m_1(z + w_2)   |   \le N^{-1/q}.
\end{equation}
To bound the second term in the sum, we used $\Im (z + \tilde w_{\overline{Y}} ) \wedge \Im (z + w_2) \ge c$ and the fact that $ |\partial_z m_1(z)| \le C\eta^{-2}$. 
This completes the proof.
\end{proof}
\bec\label{c:initialdata} Let $I$ and $\mathcal D_I$ be as in Theorem \ref{t:bulklaw}. There exist constants $p>0$ and $N_0$ such that, with overwhelming probability,
\begin{equation}
\sup_{z\in \mathcal D_I} \left|  \frac{1}{N}\tr  \left(  \frac{1}{\mathcal H(0) - z }\right)  - m_3(z)  \right| \le N^{-1/p} + \frac{N^\nu}{\sqrt{N\eta}}
\end{equation}
for $N \ge N_0$.
\eec
\begin{proof} The claim follows from Theorem \ref{t:bulklaw}, the identity $\tilde m (z) =  m_{\overline{Y}} ( z + \tilde w_X)$ from \eqref{sce3}, and the second inequality in \eqref{e:firstineq}. Since $\overline{Y}$ is diagonal, $m_{\overline{Y}} ( z + \tilde w_X)$ is sum of terms identical to the fractions in Theorem \ref{t:bulklaw}, because the latter are
\begin{equation}
\frac{ z + \tilde w_X }{  ( y_i  + (\tau-t) \widehat{A}_{ii} )^2- ( z + \tilde w_X)^2}  = \frac{1}{2} \left( \frac{1}{- z - \tilde w_X + \overline{Y}_i} + \frac{1}{- z - \tilde w_X - \overline{Y_i}} \right).
\end{equation}
\end{proof}

\section{Well posedness of dynamics}\label{s:wp}

To show the well posedness of \eqref{e:mhatsv}, it is important to ensure that the drift term, which depends on the the inverses of the eigenvalue spacings, does not become too singular. We guarantee this by adding a small perturbation to the diagonal matrix $X$ defined in \eqref{e:model}. Let
\beq
X ' = \diag(x_1, \dots, x_N) + e^{-N} Q,
\eeq
where $Q$ is an a $N \times N$ matrix of i.i.d.\ standard complex Gaussians. We first note that because the perturbation is exponentially small, it does not affect our desired conclusion. The proof is trivial and hence omitted. For the rest of this work, we use the redefined version of $\widehat M$ with $X'$ and may not explicitly indicate this.
\bel\label{l:redefinition}
If Theorem \ref{t:main} holds when $X'$ replaces $X$ in definition \eqref{e:model}, then Theorem \ref{t:main} holds. 
\eel 

We now prove the desired eigenvalue repulsion estimates. The proof of the following lemma is similar to \cite[Proposition 2.3]{che2017local}. For completeness we provide some details in the current context. 

\bel Let $P$ be a $N\times N$ matrix of complex numbers, and let $Q$ be a $N \times N$ matrix of i.i.d.\ standard complex Gaussians. Define the $2N \times 2N$ matrix $\widetilde P$ by 
\begin{equation}
\widetilde P = \bma 0 & P + e^{-N} Q \\  P^* + e^{-N} Q^* & 0  \ema
\end{equation}
Let $\gamma_1 \le \dots \le \gamma_N$ be the eigenvalues of $P$ and $\alpha_1 \le \dots \le \alpha_N$ be the positive eigenvalues of $\widetilde P$. Let $\alpha_{-i} = - \alpha_i$ denote the corresponding negative eigenvalues. Then the $\alpha_i$ are almost surely distinct, and we have the following estimates for every $\delta \in (0,1)$:
\begin{equation}
\E \sum_{i\neq j} \frac{1}{ |\alpha_i - \alpha_j| } < c_N \psi(N, P), \quad \P[ \min_{i\neq j} |\alpha_i - \alpha_j| \le \delta ] \le c_N \psi(N, P ) \delta^2,
\end{equation}
where $c_N$ is an $N$-dependent constant and 
\begin{equation}
\psi(N,P) = \exp\left(  e^{2N} \sum_{k,l} |P_{kl}|^2  \right).
\end{equation}
Finally, we have
\begin{equation}
\P[ \max_{1 \le k \le N} |\alpha_k - \gamma_k| \ge e^{-N/2} ] \le e^{ - e^{N/2} }.
\end{equation}

\eel

\begin{proof} Recall that $P$ has a singular value decomposition $P = USV^*$, where $S$ is diagonal and $U$ and $V$ are unitary. Therefore, after conjugating by the unitary block matrix  $ \frac{1}{\sqrt{2}}\bma  U & -U  \\ V & V  \ema $, which leaves invariant  the eigenvalues and the distribution of $Q$, we may suppose $P$ is  real and diagonal.

Define the index set corresponding to the off-diagonal blocks by 
\begin{equation}
\mathcal J = \{  (i,j) \colon 1 \le i,j \le 2N, \, i \le N < j \text{ or } j \le N < i   \}.
\end{equation}
Let $\mathcal H_{N}$ be set of $2N \times 2N$ Hermitian matrices with zeros in the indices $\mathcal J^c$ (the diagonal $N\times N$ blocks). We parameterize $\mathcal H_N$ by the coordinates $(w_{ij}) \in \R^{2N \times 2N}$, where $w_{ij} = 0$ if $(i,j) \in \mathcal J^c$ and $h_{ij} = w_{ij } + \ii w_{ji}$ for $j > i$ otherwise. This space is naturally equipped with the Lebesgue measure for $\R^{2 N^2}$. 

Set $\sigma_N = e^{-N}$ and write the density for $\widetilde P$ as
\begin{equation}
p_{\widetilde P}(w) = \frac{1}{Z_N} \exp \left(  - \frac{1}{2 \sigma_N^2} \sum_{(i,j) \in \mathcal J }  | w_{ij}  - \delta_{ j, i+N}  P_{i i }  |^2 \right),
\end{equation}
where we use that $P$ is real, so only the $w_{ij}$ representing the real parts of the diagonals of the off-diagonal blocks are shifted. Note the normalization constant $Z_N$ does not depend on $P$. 

In the eigenvalue--eigenvector coordinates,\footnote{The technical details of this reparameterization are similar to \cite[Proposition 2.3]{che2017local} and therefore omitted.} we have
\begin{equation}
p_{\widetilde P}(\lambda, u ,v) = \frac{1}{Z_N} \exp \left(  - \frac{1}{2 \sigma_N^2} \left( 2 \sum_{k=1}^{N} \lambda^2_{k}  + \sum_{k=1}^{N} P^2_{k k }  - 2  \sum_{k,\ell \le N }  \lambda_{\ell} \Re(u_{k\ell} v^*_{ \ell k} )P_{k k }   \right) \right) \prod_{i \neq j} ( \lambda_i - \lambda_j)^2 g(u,v),
\end{equation}
where we used 
\begin{equation}
w_{k,k+N} = \sum_{\ell=1}^N \lambda_\ell \Re(u_{k\ell} v^*_{ \ell k} )
\end{equation}
from the singular value decomposition for the upper-right block of the $(w_{ij})$ matrix. Here $g(u,v)$ is an integrable function on the compact subdomain of $\C^{N(N-1/2)} \times \C^{N(N-1/2)}$ where the map $(u,v)\rightarrow (U(u), V(v))$ taking the strictly upper triangular part of a matrix to the full Hermitian matrix is well-defined.

Using the trivial bound of $1$ on the eigenvector entries and the AM--GM inequality, we obtain
\begin{equation}
2  \sum_{k,\ell \le N }  \lambda_{\ell} \Re(u_{k\ell} v^*_{ \ell k} )P_{k k }  \le  \sum_{k=1}^N \lambda_k^2  + \sum_{k=1}^N P_{kk}^2.
\end{equation}
This implies
\begin{equation}
p_{\widetilde P}(\lambda, u ,v) \le \frac{1}{Z_N} \exp \left(  - \frac{1}{2 \sigma_N^2}  \sum_{k=1}^{N} \lambda^2_{k}    \right) \prod_{i \neq j} ( \lambda_i - \lambda_j)^2 g(u,v).
\end{equation}

Then integrating out the $g(u,v)$ term and integrating again to compute $\E \sum_{i\neq j} | \alpha_i - \alpha_j|^{-1}$, we obtain the first bound (where we use \eqref{e:hs1} and \eqref{e:hs2} to simplify the sum of the eigenvalues squared). The final inequality follows as in \cite[Proposition 2.3]{che2017local}.
\end{proof}
 
With this estimate, the following well posedness theorem is proved nearly identically to \cite[Theorem 5.2]{che2017local}. For any $\ft>0$ we define the filtration
\begin{equation}
(\mathcal F_t)_{ 0 \le t \le \ft} = ( \sigma( \widehat M(0), (B_s)_{0\le s \le t}) )_{0 \le t \le \ft},
\end{equation}
where $B_s$ is the multi-dimensional Brownian motion driving \eqref{e:mhatsv}.
 
 \bet\label{t:wp} For any $\ft > 0$, the singular values $\lambda_i(t)$ of $\widehat M(t)$ and their negatives ${- \lambda_i(t) = \lambda_{-i}(t)}$ are the unique strong solution to the equation \eqref{e:mhatsv} on $[0, \ft]$ such that 
 \begin{itemize}
 \item $\lambda(t)$ is adapted to the filtration $(\mathcal F_t)_{ 0 \le t \le \ft}$, and
 \item$ \P \left[ \lambda_{-N}(t) < \dots < \lambda_{-1}(t) < \lambda_1 (t) < \lambda_2(t) < \dots < \lambda_N(t), \text{ for almost all }t \in [ 0, \ft ] \right] =1$.
 \end{itemize}
\eet

\section{Analysis of SDEs} \label{s:sdeanalysis}

The system of SDEs for the evolution of the singular values of $\widehat{M}$ is

\begin{equation}
d \lambda_i  = \frac{1}{\sqrt{2N}} dB_i   +  \frac{1}{2 N }\sum_{j \neq i } \frac{1 - \gamma_{ij} }{\lambda_i - \lambda_j }\, dt   + R_i,
\end{equation}
for $1 \le |i| \le N$, where 
\begin{equation}
R_i =  \Re\left\langle j_i, \left(U^*(t)\widehat AV(t)-U^*(0)\widehat AV(0)\right) k_i \right\rangle\, \dd t +  \frac{1}{\sqrt{N}} \Re\left\langle j_i , \left(U^*(t)\left(\one_{(i,j)\in\mathcal{I}_\mathfrak{a}}\dd \widetilde B_{ij}\right)V(t)\right) k_i \right\rangle.
\end{equation} 
for $i \ge 1$.  We recall that with $\lambda_i$ and $\lambda_{-i}$ are coupled as discussed above so that $\lambda_i = - \lambda_{-i}$ (and the remainder terms and the $\gamma_{ij}$ are coupled in the same way). We use the redefinition noted in Lemma \ref{l:redefinition}, so that our well posedness result Theorem \ref{t:wp} applies.

Our plan is to study this system for times $0 \le t \le \tau$ with $\tau = N^{- 1 + \fb}$ and compare it to the process defined by 
\begin{equation}
d \mu_i  = \frac{1}{\sqrt{2N}} dB_i   +  \frac{1}{2 N }\sum_{j \neq i } \frac{1 }{\mu_i - \mu_j }\, dt,\qquad \mu_i(0) = \lambda_i(0),
\end{equation}
which we treat using the methods of \cite{che2019universality}. We follow closely the strategy in \cite{che2017local}, commenting on the minor differences in the current setting.

\subsection{Interpolating process}

For $0 \le \alpha \le 1$ we define the interpolating process $z_i(t,\alpha)$ by the SDE 
\begin{equation}\label{e:interp}
dz_i(t,\alpha)  = \frac{1}{\sqrt{2N}} dB_i   +  \frac{1}{2 N }\sum_{j \neq i } \frac{1 - \alpha \hat \gamma_{ij} }{z_i(t,\alpha) - z_j(t,\alpha) }\, dt   ,\quad z_i(0,\alpha) = \lambda_i(0),
\end{equation} 
with
\begin{equation}
\hat \gamma_{ij} = \gamma_{ij} \wedge N^{-\fc}.
\end{equation}
The well posedness of \eqref{e:interp} follows from the same method used to prove the well posedness of Dyson Brownian motion; see for example \cite[Lemma 4.3.3]{AGZ}.


Define
\begin{equation}
m_0(z) = \frac{1}{2N} \sum_{1 \le |i | \le N} \frac{1}{\lambda_i (0) -z},
\end{equation}
and let $m_t(z)$ be the free convolution of $m_0$ with the semicircle law at time $t$ (see \cite{schnelli2} for details):
\begin{equation}
m_t(z) = m_0(z + t m_t(z)), \quad \lim_{|z| \rightarrow \infty} m_t(z) = 0.
\end{equation}

Let $I = [-c, c]$ be the interval from Theorem \ref{t:bulklaw} on which $\mu_1^{\mathrm{sym}} \boxplus \mu_2^{\mathrm{sym}}$ has a positive density bounded above and away from zero. Let $\gamma_i^{\boxplus}$ be the $i$-th classical eigenvalue location (the $i$-th $N$-quantile) for the measure $\mu_1^{\mathrm{sym}} \boxplus \mu_2^{\mathrm{sym}}$, and define the index set $\mathcal J$ by 
\begin{equation}
\mathcal J = \{  i \colon \gamma_i^{\boxplus} \in I \}.
\end{equation}
From Corollary \ref{c:initialdata} we can deduce by standard arguments (cf. \cite[Chapter 11]{EYbook}) that there exists $c>0$ such that 
\begin{equation}
|  \gamma_i^{\boxplus} - \gamma_i(0)  | \le N^{-c}
\end{equation}
for $i\in \mathcal J$ with overwhelming probability.

The function $m_t(z)$ is the Stieltjes transform of some probability density $\rho_t(E)$. Let the classical eigenvalue locations of the free convolution $\rho_t$ be $\{ \gamma_i(t)\}_{|i| =1}^N$. Note that by the same reasoning given in \cite[Section 4.4]{che2017local}, that for any $\nu >0$ and $i,j \in \mathcal J$ with $|i-j| \ge N^\nu$,
\begin{equation}
c \frac{|i - j|} {N} \le | \gamma_i(t) - \gamma_j(t) | \le C \frac{|i - j| }{N}.
\end{equation}

The following rigidity lemmas hold. They are straightforward adaptations of the proofs of Theorem 3.1 and Corollary 3.2 of \cite{huang2016local}, and the discussion in \cite[Section 4.5]{che2017local}. The main difference is that our Brownian motions are coupled in pairs, $B_i =  - B_{-i}$. However, this does not affect the bound on the Brownian motion terms in equation (3.33) of \cite{huang2016local} in the proof the deformed law, so the same method applies here. Observe our global eigenvector bounds from Corollary \ref{e:57} are used to prove the second lemma.
\bel For any $\nu \ge 0$, it holds with overwhelming probability that
\begin{equation}
\sup_{0 \le t \le \tau } \sup_{i \in \mathcal J } | z_i(t,0) - \gamma_i(t) | \le \frac{N^\nu}{N}.
\end{equation}
\eel

\bel
With overwhelming probability for $i\in \mathcal J$,
\begin{equation}
| z_i(t,\alpha) - \gamma_i(t) | \le \frac{N^{5\fa} }{N}.
\end{equation}
\eel

\subsection{Conclusion}
The remaining stochastic analysis, including a short-range approximation and use of the energy method, is virtually identical to the argument given in \cite{che2017local}, and we obtain the following coupling.

\bep\label{p:sdecompare} Fix $\kappa >0$. Suppose that $\fb < \fa/100$ and $\fa < \fc /10$. For every time $t$ such that $0\le t \le \tau$, we have with overwhelming probability for every index $i\in \mathcal J$ that 
\beq
| \lambda_i(t) - \mu_i(t) | \le \frac{1}{N} \left(   N^{-\fc/5}   + N^{-5\fb} + N^{-1/4} \right).
\eeq
\eep

The following proposition is essentially \cite[Theorem 3.2]{che2019universality}.\footnote{It is likely that the techniques in the recent work \cite{bourgade2018extreme} could be used provide a shorter proof than the one given in \cite{che2019universality}, but since the result is already established we do not take this up here.} Compared to that reference, a certain repulsion term is present in the dynamics we study here (cf. Appendix \ref{a:realcase}), but the proof is nearly identical (and in fact strictly easier) in our case. 

We first recall the setup from that reference. Fix $\delta_1>0$ and let $g$ and $G$ be $N$-dependent parameters such that 

\beq N^{-1+\delta_1} \le g \le N^{-\delta_1},\quad G\le N^{-\delta_1}.\eeq

Let $V$ be a deterministic matrix and let $B_t=\{B_{ij}(t)\}_{1\leq i, j \leq N}$ be a matrix of i.i.d.\ standard complex Brownian motions. Define
\beq M_t = V + \frac{1}{\sqrt N } B_t,\quad H_t =  \bma  0 & M_t \\ M_t^\trans & 0 \ema . \eeq 
Let $\{s_i(t) \}_{i = -N}^N$ (omitting the zero index) be the eigenvalues of $H_t$. We set
\beq m_V(z) = \frac{1}{2N}\sum_{i= - N}^N \frac{1}{s_i(0) - z},\eeq 
where again $i=0$ is omitted in the sum.

For the next definition, we recall that $m_3(z)$ was defined as the Stieltjes transform of $\mu_1^{\mathrm{sym}} \boxplus \mu_2^{\mathrm{sym}}$.
\bed\label{def:regular}
With $g$ and $G$ as above, we say $V$ is $(g,G)$-regular with respect to $m_3$ if there exists $c>0$ such that
\beq \label{e:newreg}
\left| \Im m_V(E+i\eta)  - m_3(z)  \right| \le N^{-c} \eeq
for $z=E +  i \eta$ with  $|E| \le G$ and $\eta \in [ g , 10]$, for large enough $N$, and if there exists a constant $C_V$ such that $|v_i| \le N^{C_V}$ for all $v_i$. \eed

Let $W$ be a random matrix whose entries are i.i.d.\ complex normal variables of variance $N^{-1}$, and let $\tilde B_t=\{\tilde B_{ij}(t)\}_{1\leq i, j\leq N}$ be a matrix of i.i.d.\ standard complex Brownian motions. Define $W_t = W + N^{-1/2} \tilde{B}_t$.  Recall  $\{s_i(t)\}_{i=1}^N$ are the singular values of $M_t$, and let $\{r_i(t)\}_{i=1}^N$ be the singular values of $W_t$.

\bep\label{p:homog} Fix $\sigma>0$, and let $V$ be a deterministic matrix that is $(g,G)$-regular with respect to $m_3$. Let $M_t$, $W_t$, $\{s_i(t)\}$, and $\{r_i(t)\}$ be defined as above. Then there exists a coupling of the processes $\{s_i(t)\}$ and $\{r_i(t)\}$ such that the following holds. Given parameters $0< \omega_1 < \omega_0$ and times $t_0 = N^{-1+\omega_0}$, $t_1 = N^{-1+\omega_1}$, with the restrictions that 
\beq gN^{\sigma} \le t_0 \le N^{-\sigma} G^2,\quad 2\omega_1 < \omega_0,\eeq
there exist $C, \omega, \delta>0$  such that
\beq\label{e:desired} |s_i(t_a) - r_i(t_a) | <  C N^{-1-\delta}\eeq
with overwhelming probability for $i < N^\omega$ and $t_a=t_0+t_1$. Here $C, \omega, \delta$ are constants that depend only on $\delta_1$, $\sigma$, $\omega_0$, $\omega_1$, and the constants used to verify Definition \ref{def:regular} for $V$.
\eep

\begin{remark}
Unfortunately, the statement of \cite[Theorem 3.2]{che2019universality} omits a necessary hypothesis used in its proof. We have corrected this in the statement of Proposition \ref{p:homog}, and we now explain the changes in detail.

The definition of  $(g,G)$-regularity used in \cite[Theorem 3.2]{che2019universality} is weaker than Definition \ref{def:regular}, and merely requires that $\Im m_V$ be bounded above and below. Here, we impose the stronger condition \eqref{e:newreg}. Together with assumption \eqref{e:fix0}, \eqref{e:newreg} ensures that that $\hat \rho_{t_0}(0)$, the value at $0$ of the density corresponding to the Stieltjes transform of the free convolution of the data $V$ with $t_0$ times the semicircle law, is close to $1/ \pi$, the value at $0$ of the density of the semicircle law. The latter law governs the density of the singular values of the reference Gaussian ensemble, and this matching of densities is necessary to place the particles $V$ on the same scale as that ensemble and permit the coupling at time $t_0$ between the $s_i(t)$ and $r_i(t)$ used in the proof of \cite[Theorem 3.2]{che2019universality}. This condition on $\hat \rho_{t_0}(0)$ is tacitly assumed in the proof of \cite[Theorem 3.2]{che2019universality} but missing from its statement.

To prove Proposition \ref{p:homog}, one may follow the proof of \cite[Theorem 3.2]{che2019universality} to obtain \eqref{e:desired} up to a scaling of the particles $s_i$ by $ \pi \hat \rho_{t_0}(0)$. Using \eqref{e:newreg}, it can be shown that this scaling is $1 +  O(N^{-c})$, and we obtain the claimed result. The details of the latter argument can be found in the proof of \cite[Theorem 2.4]{che2017local}; see the discussion starting above (4.110).

\end{remark}

The hypotheses of Proposition \ref{p:homog} are verified with overwhelming probability for the singular values of $\widehat M(0)$ by Corollary \ref{c:initialdata}. Combining Proposition \ref{p:homog} with Proposition \ref{p:sdecompare}, we obtain short-time relaxation of the singular value dynamics. 

\bet\label{t:shorttime} Fix $\sigma>0$, $\kappa>0$, suppose that $\fb < \fa/100$ and $\fa < \fc /10$, and retain the definitions of Proposition \ref{p:sdecompare}.  Then there exists a coupling of the processes $\{\lambda_i(t)\}$ and $\{r_i(t)\}$ and a constant $N_0(\sigma, \kappa, \fa, \fb, \fc)$ such that 
\beq |\lambda_i(t_a) - r_i(t_a) | < N^{-1-\delta}\eeq
with overwhelming probability for $i < N^\omega$ and $N \ge N_0$.
\eet
We are now positioned to prove our main theorem.\\

\noindent {\bf Proof of Theorem \ref{t:main}.}  Setting $t_a = \tau$ in Theorem \ref{t:shorttime}, we have $|\lambda_1(\tau) - r_1(\tau) | < N^{-1-\delta}$ where $\lambda_1$ is the least singular value of $\widetilde M(\tau) = M(\tau)$ and $r_1(\tau)$ is the least singular value of a matrix with distribution $\sqrt{1 + \tau} W$, where $W$ is a matrix of i.i.d.\ standard complex Gaussians. Note that $M(\tau)$ has the same law as $M$, so if $\lambda_1(M_N)$ is the least singular value of $M_N$, we have 
\beq\label{e:1coupling} |\lambda_1(M_N) - r_1(\tau) | < N^{-1-\delta}.\eeq
The distribution of the least singular value of a Gaussian matrix is known explicitly. For $W$  and any $r \ge 0$ \cite{E88}, 
\beq \label{e:exactdist} \P( N \lambda_1(W)  \le r ) = 1 - e^{- r^2}.
\eeq
Therefore, 
\beq \P( N \lambda_1(\sqrt{1 + \tau}W)  \le \sqrt{1 + \tau} r ) = 1 - e^{- r^2}.\eeq
We now show the $\sqrt{1 + \tau}$ factor is negligible, so that we may compare $\lambda_1(M_N)$ directly to $\lambda_1(W)$. We compute, using $1 - e^{-x} \le x$,
\begin{multline}  \left| \P( N \lambda_1(W)  \le r ) - \P( N \lambda_1(\sqrt{1 + \tau}W)  \le r ) \right| = \left| e^{- r^2/(1+\tau)}   - e^{- r^2} \right| \\ = e^{-r^2/(1+\tau)} \left| 1 - \exp\left(   - r^2 \left( \frac{\tau}{1+\tau}  \right) \right)   \right| \le \left[ e^{ - \frac{r^2}{1+\tau} }  \frac{r^2}{1+\tau } \right] \cdot \tau \le C \tau = O(N^{-c}).
\end{multline}

By \eqref{e:1coupling}, 
\beq
 \P( N r_1(\tau)  \le  r - N^{-\delta} ) \le  \P ( N \lambda_1(M_N) \le r) \le \P( N r_1(\tau)  \le  r + N^{-\delta} ) 
\eeq
We deduce 
\beq  \P( N \lambda_1(W)  \le  r - N^{-\delta} )  - CN^{-c} \le  \P ( N \lambda_1(M_N) \le r) \le \P( N \lambda_1(W)  \le  r + N^{-\delta} )  + C N^{-c} .\eeq
By \eqref{e:exactdist}, $N\lambda_1(W)$ has a bounded density, so the $N^{-\delta}$ terms in the above may be removed with $O(N^{-c})$ error, and we conclude that 
\beq
 | \P ( N \lambda_1(M_N) \le r) -  \P ( N \lambda_1(W ) \le r)| = O(N^{-c})
\eeq
as desired. \ep
 
\appendix

\section{Derivation of dynamics}  \label{a:ito}
The following is a formal calculation that ignores the technical issue of possible eigenvalue collisions. It is used in Section \ref{s:wp}, where this issue is dealt with rigorously. 
\subsection{Calculation} With $\widehat M$ as above, we define the $2N \times 2N$ block matrix
\begin{equation}\label{e:mx}
X = \bma 0 & \widehat M \\ \widehat M^\trans & 0  \ema =  \bma 0 & M(t) + (\tau-t)U^*(0)\widehat AV(0)  \\ M(t)^* + (\tau-t)V(0)^*\widehat AU(0) & 0  \ema.
\end{equation}
Observe that the eigenvalues of $X$ are the singular values of $\widehat M$ and their negatives. Let $\widehat M =  J  S  K^*$ be the singular value decomposition of $\widehat M$. Then a matrix of normalized eigenvectors for $X$ is 
\begin{equation}\label{e:symX}
H =  \frac{1}{\sqrt{2} } \bma  J & - J \\ K & K \ema.
\end{equation}
We follow the approach of \cite[Chapter 12]{EYbook} to compute the dynamics of the eigenvalues of $X$. Denote the eigenvalues of $X$ by $\lambda_\alpha$ with corresponding eigenvectors $u_\alpha$. For the elements $x_{ij}$ of $X$ that are not identically zero, we have 
\begin{equation}
\frac{\partial \lambda_\alpha}{\partial x_{ij}} = u^*_\alpha (i) u_\alpha(j),\quad \frac{\partial u_\alpha(i)}{ \partial x_{k l} } = \sum_{\beta \neq \alpha} \frac{u^*_\beta(k) u_\alpha(l)}{\lambda_\alpha - \lambda_\beta} u_\beta(i),
\end{equation} 
and by the chain rule,
\begin{equation}
\frac{\partial^2 \lambda_\alpha}{\partial x_{kl} \partial x_{ij} } = \sum_{\beta \neq \alpha} \frac{1}{\lambda_\alpha - \lambda_\beta} \left[ u^*_\alpha(k) u_\beta(l) u^*_\beta(i) u_\alpha(j)  +  u^*_\beta(k) u_\alpha(l) u^*_\alpha(i) u_\beta(j)  \right].
\end{equation}
It\^o's formula gives
\begin{equation}
d\lambda_\alpha = \sum_{i,j}  \frac{\partial \lambda_\alpha}{\partial  x_{ij} }  d x_{ij} + \frac{1}{2} \sum_{i,j,k,l} \frac{\partial^2 \lambda_\alpha}{\partial x_{kl} \partial x_{ij} }  (dx_{ij}) (dx_{kl}).
\end{equation}
The first term is, for $\alpha \le N$, 
\begin{align}
 \sum_{i,j}  \frac{\partial \lambda_\alpha}{\partial  x_{ij} }  d x_{ij} &= \sum_{i,j} u^*_\alpha(i)  u_\alpha(j) dx_{ij}  =  u_\alpha^* (dX) u_\alpha \\
 &= \frac{1}{2 } j_\alpha^* \left( \frac{1}{ \sqrt{N}} d \widehat B   + (U^* \widehat AV-U^*(0)\widehat AV(0))\, \dd t - \frac{1}{\sqrt{N}}U^*(\one_{(i,j)\in\mathcal{I}^c_\mathfrak{a}}\, \dd \widetilde B_{ij})V \right) k_\alpha\\
 &+ \frac{1}{2 }  k_\alpha^* \left( \frac{1}{ \sqrt{N}} d \widehat B   + (U^*\widehat AV-U(0)\widehat AV(0)^*)\, \dd t - \frac{1}{\sqrt{N}}U^*(\one_{(i,j)\in\mathcal{I}^c_\mathfrak{a}}\, \dd \widetilde B_{ij})V \right)^* j_\alpha.
\end{align}
We see that 
\begin{equation}
 \frac{1}{2  \sqrt{N}} \left(  j_\alpha^* d \widehat B k_\alpha +  k_\alpha^* d \widehat B^* j_\alpha \right) =  \frac{1}{\sqrt{2 N}} d B_\alpha,
 \end{equation}
 where $\{dB_\alpha\}_{\alpha=1}^N $ is a set of independent standard real Brownian motions. The independence follows from an explicit computation, noting that $(d\widehat B_{ij})( d\widehat B_{kl} )= \delta_{il} \delta_{jk}$. The remaining terms are
\beq   \Re\langle j_\alpha, (U\widehat AV^*-U(0)\widehat AV(0)^*) k_\alpha \rangle\, \dd t +  \frac{1}{\sqrt{N}} \Re\langle j_\alpha, (U(\one_{(i,j)\in\mathcal{I}_\mathfrak{a}}\dd \widetilde B_{ij})V^*) k_\alpha \rangle . \eeq
The second term is 
\begin{multline}
\frac{1}{2} \sum_{i,j,k,l} \frac{\partial^2 \lambda_\alpha}{\partial x_{kl} \partial x_{ij} }  (dx_{ij}) (dx_{kl}) \\= \frac{1}{2} \sum_{i,j,k,l} \sum_{\beta \neq \alpha} \frac{1}{\lambda_\alpha - \lambda_\beta} \left[ u^*_\alpha(k) u_\beta(l) u^*_\beta(i) u_\alpha(j)  +  u^*_\beta(k) u_\alpha(l) u^*_\alpha(i) u_\beta(j)  \right] \, dx_{ij}  dx_{kl}.
\end{multline}
The first contribution is 
\begin{multline}
 \frac{1}{2} \sum_{i,j,k,l}\sum_{\beta \neq \alpha} \frac{1}{\lambda_\alpha - \lambda_\beta} \left[ u^*_\alpha(k) u_\beta(l) u^*_\beta(i) u_\alpha(j)  +  u^*_\beta(k) u_\alpha(l) u^*_\alpha(i) u_\beta(j)  \right]   \frac{1}{ N } (d \widehat B'_{ij}  )( d \widehat B'_{kl}),
\end{multline}
where 
\begin{equation}
d \widehat B'  = \bma  0 & d\widehat B \\ d\widehat B^* & 0  \ema .
\end{equation}
This vanishes unless $i=l$, $j=k$, and exactly one of $i$ or $j$ is greater than $N$, due to the covariation factor. Summing over $i$ and $j$, we obtain the norm of the first or last half of each $u_\alpha$, that is  $\| j _\alpha \|_2 /2$ or $\| k _\alpha \|_2 /2$, both of which are $1/2$. We then recover the drift term
\begin{equation}
\frac{1}{2N} \sum_{\beta \neq \alpha} \frac{dt}{\lambda_\alpha - \lambda_\beta}.
\end{equation}
The remaining contribution is 
\begin{multline}\label{e:remaining}
 - \frac{1}{2 N } \sum_{i,j,k,l}\sum_{\beta \neq \alpha} \frac{1}{\lambda_\alpha - \lambda_\beta} \left[ (u^*_\alpha(k) R'_{kl} u_\beta(l))( u^*_\beta(i) R' _{ij} u_\alpha(j) ) +  (u^*_\beta(k) R'_{kl} u_\alpha(l) )(u^*_\alpha(i) R' _{ij} u_\beta(j) ) \right] .
\end{multline}
where 
\begin{equation}
R' = \bma 0 &  U^*(\one_{(i,j)\in\mathcal{I}^c_\mathfrak{a}}\, \dd \widetilde B_{ij})V \\ ( U^*(\one_{(i,j)\in\mathcal{I}^c_\mathfrak{a}}\, \dd \widetilde B_{ij})V)^* & 0 \ema.\end{equation}
We perform the sum on $i$ and $j$ first. We have
\begin{equation}\label{e:factors}
\sum_{i,j} u^*_\alpha(i) R' _{ij} u_\beta(j) = \frac{1}{2} \left( j_\alpha^* U^*(\one_{(i,j)\in\mathcal{I}^c_\mathfrak{a}}\, \dd \widetilde B_{ij})V k_\beta +  k_\alpha^* V^*(\one_{(i,j)\in\mathcal{I}^c_\mathfrak{a}}\, \dd \widetilde B_{ij})^* U j_\beta\right).
\end{equation}
Define the  column vectors
\begin{equation}
w_\alpha = U j_\alpha , \quad z_\alpha = V  k_\alpha,
\end{equation}
and set $R = (\one_{(i,j)\in\mathcal{I}^c_\mathfrak{a}}\, \dd \widetilde B_{ij})$. Then since the quadratic variation of a standard complex Brownian motion is zero, and the elements of $R$ are independent,
\begin{align}
(w^*_\alpha R z_\beta  + z^*_\alpha R^* w_\beta  )(w^*_\beta R z_\alpha  + z^*_\beta R^* w_\alpha ) &= w^*_\alpha R z_\beta z_\beta^* R^* w_\alpha + z^*_\alpha R^* w_\beta w_\beta^*  R z_\alpha\\ &= \sum_{(i,j) \in \mathcal{I}_\mathfrak{a}^c} |w_\alpha(i)|^2 |z_\beta(j)|^2    +  \sum_{(i,j) \in \mathcal{I}_\mathfrak{a}^c}  |w_\beta(i)|^2 |z_\alpha(j)|^2  \\ & \defeq 2 \gamma_{\alpha \beta}.
\end{align} 
Then \eqref{e:remaining} becomes, remembering the factors of $1/2$ from \eqref{e:factors} and using $\gamma_{\alpha\beta} = \gamma_{\beta\alpha}$, 
\begin{equation}
 -  \frac{1}{2 N }\sum_{\beta \neq \alpha} \frac{\gamma_{\alpha \beta} \, dt}{\lambda_\alpha - \lambda_\beta} .
\end{equation}

We obtain the following SDE for the eigenvalues of $X$, which are the singular values of $\widehat M$ and their negatives. We label the positive eigenvalues by $\{\lambda_i \}_{i = 1}^N$ and the negative eigenvalues by  $\{\lambda_i \}_{i  = - 1}^{-N}$, where we have set $\lambda_{-i} = - \lambda_{i}$. The same convention holds for the Brownian motions $B_i$: there are $2N$ of them, and the ones with positive indices are coupled to those with negative indices by $B_{- i } = - B_{ i }$. The final SDE is, for $i >0$:
\begin{equation}
d \lambda_i  = \frac{1}{\sqrt{2N}} dB_i   +  \frac{1}{2 N }\sum_{j \neq i } \frac{1 - \gamma_{ij} }{\lambda_i - \lambda_j }\, dt   + R_i,
\end{equation}
where 
\begin{equation}
R_i =  \Re\langle j_i, (U^*\widehat AV-U^*(0)\widehat AV(0)) k_i \rangle\, \dd t +  \frac{1}{\sqrt{N}} \Re\langle j_i , (U^*(\one_{(i,j)\in\mathcal{I}_\mathfrak{a}}\dd \widetilde B_{ij})V) k_i \rangle.
\end{equation}
For $i < 0$ one can check that $R_{i} = - R_{-i}$ and $\gamma_{ij} = - \gamma_{-i,j}$.

In preparation for the next section, we note that when we conjugate the initial data by orthogonal instead of unitary matrices, the matrix \eqref{e:mx} is a function of $N^2$ real variables instead of $2N^2$ real variables ($N^2$ complex entries). In this case we obtain (see, for example, \cite[Chapter 12]{EYbook}) that
\begin{equation}
\frac{\partial \lambda_\alpha}{\partial x_{ij}} = 2 u_\alpha (i) u_\alpha(j),\quad \frac{\partial u_\alpha(i)}{ \partial x_{k l} } = \sum_{\beta \neq \alpha} \frac{u_\beta(k) u_\alpha(l) + u_\beta(l) u_\alpha(k)}{\lambda_\alpha - \lambda_\beta} u_\beta(i),
\end{equation} 
where we view $X$ as a function of $N^2$ real variables $x_{ij}$ with $1\le i \le N$ and $N \le j \le 2N$. Using the representation \eqref{e:symX}, we find
\begin{align}
u_{-\alpha} (k) u_\alpha(l) + u_\beta(l) u_\alpha(k) &= u_{-\alpha} (k) u_\alpha(l) + u_{-\alpha} (l) u_\alpha(k)\\ &=  - u_{\alpha} (k) u_\alpha(l) + u_{\alpha} (l) u_\alpha(k)\\ &= 0.
\end{align}
Therefore, we see that the sum in the drift component now omits the term with $\beta = - \alpha$, and there is no repulsion between $\lambda_{\alpha}$ and $\lambda_{-\alpha}$. The rest of the derivation is completed as before. 
\section{Real case}\label{a:realcase}

We now consider the real analogue of the model of Section \ref{s:model}, where the initial data is conjugated by orthogonal matrices. Precisely, in this section we consider the matrix ensemble
\beq \label{e:realmodel}
M = R^\trans X T +  U^\trans YV , 
\eeq
where $X=\diag(x_1,\dots , x_N)$ and $Y=\diag(y_1, \dots , y_N)$ are deterministic diagonal matrices and $R,T,U,V$ are independent and distributed according to the Haar measure on the orthogonal group $O(N)$. We retain the hypothesis \eqref{e:xybound} and the assumptions labeled (1) through (7) on $X$ and $Y$ given in Section \ref{s:model}. 

The least singular value in the real case displays qualitatively different behavior than its counterpart in the complex case, as indicated by the accompanying simulation results. The density for $\lambda_1$ vanishes at zero in the complex model, but remains positive in the real model. The singular value distribution in the real case is said to have a \emph{hard edge} at zero.

\begin{figure}
    \centering
    \begin{minipage}{0.5\textwidth}
        \centering
        \includegraphics[width=\textwidth]{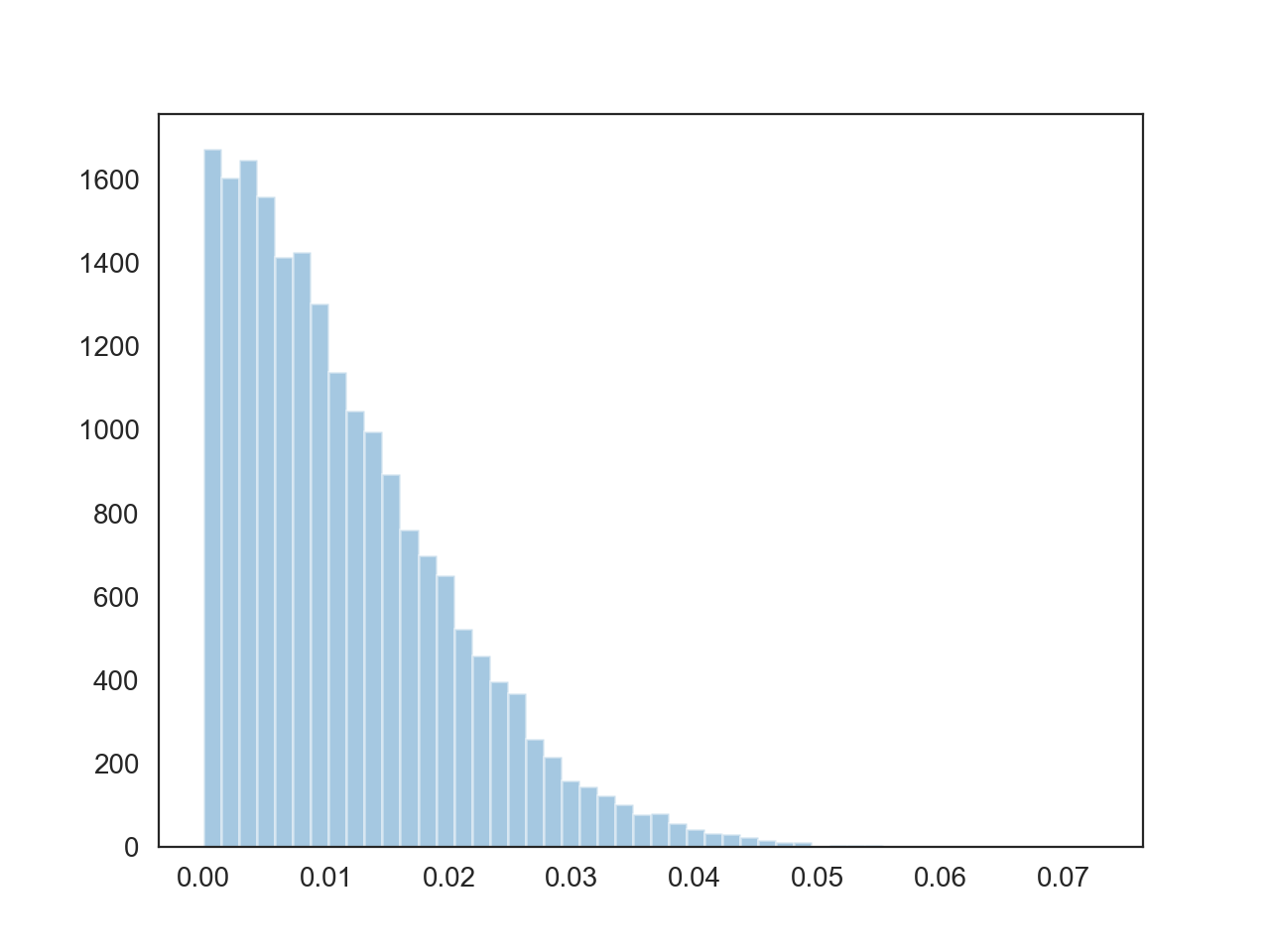} 
        \caption{Simulated distribution of the least singular value of the real model, with the elements of $X$ and $Y$ chosen uniformly from $[0,1]$, matrix size $N=200$, and $2\times 10^4$ samples.}
    \end{minipage}\hfill
    \begin{minipage}{0.5\textwidth}
        \centering
        \includegraphics[width=\textwidth]{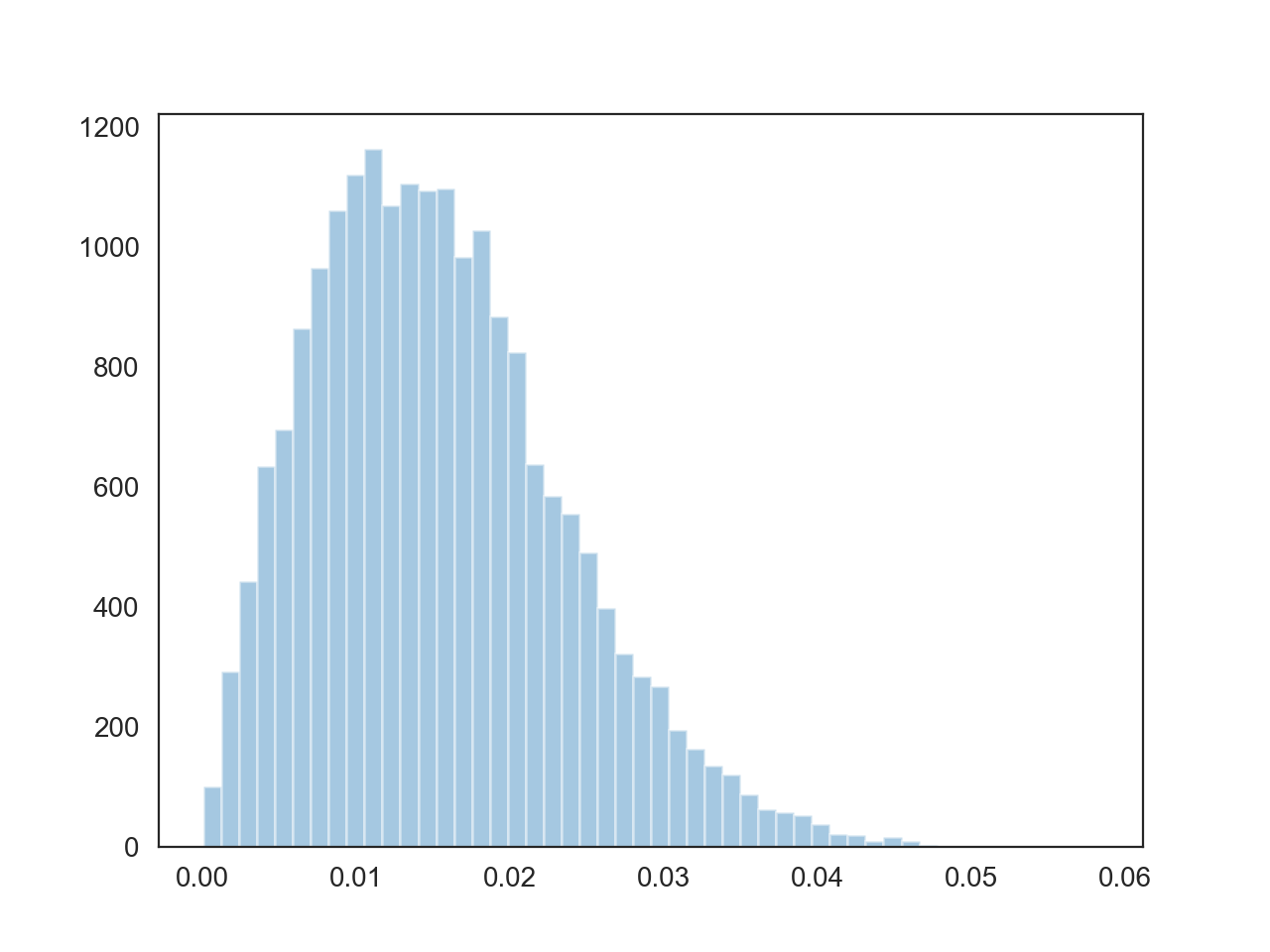} 
        \caption{Simulated distribution of the least singular value of the complex model, with the elements of $X$ and $Y$ chosen uniformly from $[0,1]$, matrix size $N=200$, and $2\times 10^4$ samples.}
    \end{minipage}
\end{figure}

This phenomenon may be understood dynamically. As discussed in Appendix \ref{a:ito}, the drift term in the complex case has the repulsion component 
 \beq \frac{1}{2 N }\sum_{j \neq i } \frac{1 - \gamma_{ij} }{\lambda_i - \lambda_j }\, dt ,\eeq
while the same computation in the real case yields the repulsion term 
  \beq \frac{1}{2 N }\sum_{j \neq i, -i } \frac{1 - \gamma_{ij} }{\lambda_i - \lambda_j }\, dt ,\eeq
  with the interaction between $\lambda_{i}$ and $\lambda_{-i}$ removed. For $\lambda_{1}$, this means there is no force from $\lambda_{-i}$ pushing it away from the origin, resulting in the hard edge.

The model \eqref{e:realmodel} can be handled by the same method used for \eqref{e:model}. The definition of the matrix dynamics in Section \ref{s:dynamicsdefinition} is the same except for obvious changes, such as the use of orthogonal matrices and real symmetric Brownian motions. This leads to virtually the same singular value dynamics as in Appendix \ref{a:ito}, with the important exception of the interaction term noted above. The estimates of Section \ref{s:estimates} are also essentially unchanged. An inspection of the proofs referenced in Section \ref{s:wp} and Section \ref{s:sdeanalysis} shows that they still apply to the dynamics in the real case. An important point is that the short-time universality result Proposition \ref{p:homog} still holds without the regularizing force from $\lambda_{-1}$; this was the original form of the result stated in \cite{che2019universality}. Finally, for the exact form of the distribution of the least singular value for the Gaussian matrix, we use the form with quantitative error given in \cite[Theorem 1.3]{leastsv}.

We obtain the following analogue of Theorem \ref{t:main} for the real model. 
\bet\label {t:main2}
Let $\lambda_1(M_N)$ be the least singular value of the random matrix ensemble \eqref{e:realmodel}. For all $r\ge 0$, we have 
\beq \P (N \lambda_1 ( M_N) \le r ) = 1  - e^{ - r^2/2 -r } + O(N^{-c})
\eeq
where $c>0$ is an absolute constant uniform in $r$.
\eet

\section{Preliminary estimates} \label{s:estimates}

\subsection{Diagonal matrices}

\bep\label{prop:empiricalbounds}
	For any $a>0$ and $E\in\R$, we have
	\begin{equation}
		\frac{1}{2N}\sum_{\absa{y_i-E}\geq N^{-1+a}} \frac{1}{\absa{y_i-E}^2} \leq 2C_a N^{1-a} ,
	\end{equation}
	\begin{equation}
		\frac{1}{2N}\sum_{\absa{y_i-E}\geq N^{-1+a}} \frac{1}{\absa{y_i-E}} \leq 2C_a\log N+ 4,
	\end{equation}
	where the sums are taken over indices $i$ such that $1\le |i| \le N$. \eep
\begin{proof}
Let $\eta=N^{-1+a}$. Note that
\begin{equation}
	\frac{1}{2N}\sum_{\absa{y_i-E}\geq \eta} \frac{\eta}{\absa{y_i-E}^2}  \leq \frac{1}{N}\Im \sum_{|i|=1}^N \frac{1}{y_i-(E+\ii\eta)} \leq 2C_a.
\end{equation}
Divide both sides by $\eta$ to obtain
	\begin{equation}
		\frac{1}{2N}\sum_{\absa{y_i-E}\geq \eta} \frac{1}{\absa{y_i-E}^2} \leq \frac{2 C_a}{\eta} = 2 C_a N^{1-a}.
	\end{equation}
This proves the first inequality in the proposition. For the second inequality, note that for $x>\eta$, we have
\begin{equation}
	\frac{1}{x}\leq \frac{2}{x+1} +  \int_\eta^1 \frac{2\, \dd t}{(x+t)^2}  \leq 2 +  \int_\eta^1\frac{2\, \dd t}{x^2+t^2} .
\end{equation}
Taking $x=\abs{y_i-E}$ and summing over $i$ such that $\absa{y_i-E}\geq \eta$, we have
\begin{equation}	
	\frac{1}{2N}\sum_{\absa{y_i-E}\geq \eta} \frac{1}{\absa{y_i-E}} \leq  4 +  \frac{1}{2N} \int_{\eta}^1\sum_{|i| =1}^N\frac{4\, \dd t}{\absa{y_i-(E+\ii t)}^2} .
\end{equation}
Using $\frac{1}{\absa{y_i-(E+\ii t )}^2} = \frac{1}{t } \Im\left( \frac{1}{y_i-(E+\ii t)}\right)$ and the hypothesized bound \eqref{e:mYbound} on $m_Y(E+\ii\eta)$, we have
\begin{equation}	
	\frac{1}{2N}\sum_{\absa{y_i-E}\geq \eta} \frac{1}{\absa{y_i-E}} \leq 4 +C \int_\eta^1 \frac{\dd t}{t} \leq 4 + C\absa{\log \eta}.
\end{equation}
Here $C=2\sup_{\eta\leq t\leq 1} \absa{m_Y(E+\ii t)}$, which is bounded by $2C_a$ according to the maximum principle for holomorphic functions. Recall that $\eta=N^{-1+a}$, so we have
\begin{equation}
	\frac{1}{2N}\sum_{\absa{y_i-E}\geq \eta} \frac{1}{\absa{y_i-E}}\leq 2C_a\log N+4.
\end{equation}\end{proof}

Recalling definition \eqref{ifa} and using Proposition \ref{prop:empiricalbounds} with $E= y_i$, we immediately have a bound for $A$:
\be\label{e:Abound}
	\norm{A} = \max_{1\leq i\leq N} \abs{A_{ii}} \leq CN^{1-\mathfrak{a}}.
\ee
Similarly we obtain
\be\label{e:hatAbound}
	\norm{\widehat A} = \max_{1\leq i\leq N} \abs{\widehat A_{ii}} \leq C( 1 + \log N).
\ee \eer

\subsection{Unitary flow} The proofs of the following lemmas are essentially identical to those of \cite[Theorems 3.1 and 3.3]{che2017local}.\footnote{We recall that equation (3.7) in this reference is derived by applying the formula $ \frac{d}{dt} e^{\theta X(t)}  = \int_0^\theta e^{\alpha X(t)} \frac{dx(t)}{dt} e^{(\theta - \alpha) X(t) }\, dt $, which holds for any one-parameter matrix subgroup $X(t)$ \cite{wilcox1967exponential},  to compute the derivative of the matrix exponential with respect to each matrix entry, in conjunction with It\^{o}'s formula.}
\bel
For $\fa,\fb,U,$ as above, 
\begin{equation}
\P [ \sup_{0 \le t \le \tau } \norm{U(t) - I } \ge N^{-10 \fb }  ] \le \exp\left( - N^{10 \fb} \right),
\end{equation}
and the same estimate holds for $V$.
\eel

For any $t_0 \le \tau$ define $\widehat U(t_0) = U(t) U(t_0)^*$.

\bel\label{l:holder} For $N$ large enough the following holds. For any $0 \le t_0 \le t \le \tau$, $| t - t_0| \le 1/N$,
\begin{equation}
\P \left[ \sup_{t_0 \le s \le t } \norm{\widehat U(s) - \widehat U(t_0) }  \ge (N(t - t_0))^{1/4}  \right] \le \exp \left(   - N^{\fa /3} \right).
\end{equation}
Also, for any $0\le t_0 \le t \le \tau$, $|t - t_0| \le r \le 1/N$, 
\begin{equation}
\P\left[  \sup_{t_0 \le s \le t }  \norm{ \widehat U(s) - \widehat U(t_0) } \ge r^{9/20}   \right] \le \exp\left(  - c_N r^{-1/10}  \right).
\end{equation}
where $c_N > 0$ depends on $N$. 
\eel

\subsection{Sufficient conditions for positive density}

The next lemma follows from the argument in \cite[Lemma 3.2]{bao2017spectral}. We provide the reasoning again here for completeness.

\bel
Let $\mu_\alpha, \mu_\beta$ be probability measures with density functions $\rho_\alpha, \rho_\beta$ that are symmetric about zero and are strictly positive on $[-r_0, r_0]$ for some $r_0>0$. Then $\mu_\alpha \boxplus \mu_\beta$ has a density, and that density is bounded above and away from zero in a neighborhood of zero.
\eel

\begin{proof}
According to \cite[Corollary 8]{Be14}, $\mu_\alpha \boxplus \mu_\beta$ has a bounded density. It remains to show it is bounded away from zero. By Proposition \ref{p:extension}, the corresponding subordination functions $w_\alpha, w_\beta$ extend continuously to $0$ with values in $\C^+ \cup \R \cup \{\infty\}$. By the equations defining the free convolution, it suffices to show these limits are not infinite to show that the density $\mu_\alpha \boxplus \mu_\beta$ is bounded below in a neighborhood of $0$. 

We proceed by contradiction. Fix $r < r_0/2 $ and define 
\begin{equation}
\mathcal E = \left\{ z \in C^+ \cup \R \colon |z| \le r  \right\}.
\end{equation}
Let $L > r_0$ and $M > 10$ be large parameters to be fixed later. We first suppose that there exists $z \in \mathcal E$ such that $|w_\alpha (z)| > LM$ and $|w_\beta(z)| > L$. The defining equations for the free convolution give
\begin{equation}
( w_\alpha + w_\beta - z )^{-1} = \int_\R \frac{ d \mu_\beta(x)}{w_\alpha - x } = w_\alpha^{-1} + O( w_\alpha^{-2}),
\end{equation}
where the $O$ notation is with respect to the limit $L \rightarrow \infty$. The above equation gives 
\begin{equation}
\frac{w_\beta}{w_\alpha} = O(w_\alpha^{-1}). 
\end{equation} 
This contradicts $L/|w_\alpha| \le |w_\beta/w_\alpha|$ (which holds by our assumptions on $w_\alpha,w_\beta$) for $L$ sufficiently large. 

We next suppose $|w_\alpha (z)| > LM$ and $|w_\beta(z)| \le  L$, and find from the definition of free convolution that for $z\in \mathcal E$ and $M$ sufficiently large, 
\begin{equation}
\frac{1}{|m_\alpha(w_\beta)|} = |w_\alpha + w_\beta -z | \ge \frac{ML}{2}.
\end{equation}
By symmetry of $\mu_\alpha$ and $\mu_\beta$ we know that $w_\beta$ is imaginary for $z$ on the imaginary line $\{i\eta|\eta\in\R\}$. But $m_\alpha(z)$ has no zeros on the imaginary line, as $\rho_\alpha$ is positive near $0$.  So it is bounded away from zero in $z \in \mathcal E$. For $M$ large we reach a contradiction. This completes the proof.
\end{proof}

In the case $\mu_\alpha = \mu_\beta$, only the first part of the previous argument is required.

\bel
Let $\mu_\alpha$ be a symmetric probability measure, not necessarily absolutely continuous, supported at more than 2 points. Then $\mu_\alpha \boxplus \mu_\alpha$ has a density, and that density is bounded above and away from zero in a neighborhood of zero.
\eel


\bibliography{addition}{}
\bibliographystyle{abbrv}

\end{document}